\documentclass[11pt,a4paper]{amsart}
\usepackage[T1]{fontenc}
\usepackage[utf8]{inputenc}
\usepackage{lmodern,microtype,mathtools,amssymb,amsmath}
\usepackage[margin=25mm]{geometry}
\usepackage{enumitem}
\usepackage{xcolor}
\usepackage[colorlinks=true,linkcolor=blue!40!black,citecolor=blue!40!black,urlcolor=blue!40!black]{hyperref}

\newtheorem{theorem}{Theorem}[section]
\newtheorem{proposition}[theorem]{Proposition}
\newtheorem{lemma}[theorem]{Lemma}
\newtheorem{corollary}[theorem]{Corollary}
\theoremstyle{definition}

\theoremstyle{remark}
\newtheorem{remark}[theorem]{Remark}
\numberwithin{equation}{section}
\newcommand{\T}{\mathbb T}
\newcommand{\R}{\mathbb R}
\newcommand{\Z}{\mathbb Z}
\newcommand{\C}{\mathbb C}
\newcommand{\av}{\mathcal A}
\newcommand{\HH}{\mathcal H}
\newcommand{\ZZ}{\mathcal Z}
\newcommand{\CC}{\mathcal C}
\newcommand{\ii}{\mathrm i}
\newcommand{\dd}{\,d}
\newcommand{\norm}[1]{\left\|#1\right\|}
\newcommand{\abs}[1]{\left\vert#1\right\vert}
\newcommand{\cform}{\mathfrak b}
\newcommand{\err}{\mathcal B}
\newcommand{\deckav}{\mathcal A_{\rm deck}}

\DeclareMathOperator{\Id}{Id}
\DeclareMathOperator{\Div}{div}
\DeclareMathOperator{\curl}{curl}
\setlist[enumerate]{label=\textnormal{(\roman*)},leftmargin=2em,itemsep=4pt}
\setlist[itemize]{leftmargin=1.7em,itemsep=3pt}
\allowdisplaybreaks[2]
\title[Grad counterexamples of any knot type]{Counterexamples to Grad's conjecture\\ of any knot type}
\author[A. Enciso]{Alberto Enciso}
\address{Instituto de Ciencias Matem\'aticas, Consejo Superior de
Investigaciones Cient\'\i ficas, 28049 Madrid, Spain}
\email{aenciso@icmat.es}
\author[J. Pe\~nafiel-Tom\'as]{Javier Pe\~nafiel-Tom\'as}
\address{Institute of Mathematics, \'Ecole polytechnique F\'ed\'erale
de Lausanne, 1015 Lausanne, Switzerland}
\email{javier.penafieltomas@epfl.ch}
\author[D. Peralta-Salas]{Daniel Peralta-Salas}
\address{Instituto de Ciencias Matem\'aticas, Consejo Superior de
Investigaciones Cient\'\i ficas, 28049 Madrid, Spain}
\email{dperalta@icmat.es}

\begin{document}
\begin{abstract}
We construct smooth magnetohydrostatic counterexamples to Grad's conjecture in toroidal domains of arbitrary knot type. The magnetic field is nonvanishing, and nested pressure tori fill the domain outside a nondegenerate magnetic axis. The equilibria have trivial Euclidean symmetry group and occur in smooth families of inequivalent configurations. On each regular pressure torus the magnetic flow is smoothly conjugate to a linear flow with any prescribed positive constant slope. A variation of this argument gives equilibria of arbitrary knot type with twist (that is, nonconstant rotational transform) and any finite regularity. A key ingredient of the proof is a deformation method that creates knot topology from axial symmetry by separating the sheets of a covering immersion into a braid, which may be useful in other deformation problems.
\end{abstract}
\maketitle

\section{Introduction}\label{kg1}

The celebrated Grad's conjecture~\cite{Grad} predicts symmetry restrictions for (non-isolated) smooth magnetic equilibria with nested pressure surfaces. More precisely, using the formulation given by Constantin, Drivas and Ginsberg~\cite{CDG}, the conjecture is the following: a smooth non-isolated equilibrium on a toroidal or cylindrical domain, nonvanishing away from its magnetic axis and foliated by pressure surfaces, should have plane-reflection, axial or helical symmetry. Here non-isolation means that nearby equilibria exist beyond those obtained by translations and rescalings.

There are several constructions of smooth equilibria outside axial symmetry. Lortz~\cite{Lortz} obtained toroidal equilibria with closed magnetic lines and a plane-reflection symmetry, but the pressure level sets fail to define a foliation by nested toroidal surfaces. More recently, Drivas, Elgindi and Ginsberg~\cite{DEG} constructed smooth families in periodic cylindrical geometry with a nondegenerate pressure foliation and no continuous Euclidean symmetry. Their equilibria retain a plane reflection and a finite rotational symmetry, and all their trajectories are closed. See also~\cite{W,SK,KS} for similar constructions.

Smooth counterexamples to Grad's conjecture in embedded solid tori have been obtained very recently. First, in a technical tour de force, G\'omez-Serrano, Liehr and Taylor~\cite{GSLT} obtained smooth families whose full Euclidean symmetry group, allowing reversal of the field, is the cyclic group~$C_N$ for every integer~$N$ sufficiently large. Their magnetic field vanishes on a circular axis and all its nonconstant trajectories are periodic, winding once around the circular axis (poloidal). Remarkably, Landreman~\cite{Landreman} obtained explicit non-axisymmetric analytic families with nested pressure tori, including one with constant rotational transform and another with nonconstant transform, in which rational and irrational surfaces coexist. Issan, Citrin, Landreman and Tracey~\cite{ICLT} extended these analytic families to remove stellarator symmetry, and Issan, Jos\'e, Landreman, Citrin and Tracey~\cite{IJLCT} constructed explicit non-axisymmetric families with arbitrary integer rotational transform (but keeping the stellarator symmetry). Analytic counterexamples with trivial symmetry group have been obtained by Lacy and Sato~\cite{LS}. All of these results concern symmetry and magnetic dynamics in unknotted toroidal domains.

In this article we are concerned with the construction of counterexamples to Grad's conjecture in toroidal domains of arbitrary knotted topology. Prescribing the knot type leads to many difficulties. Indeed, the Beltrami fields constructed in~\cite{EP} have vortex tubes of arbitrary knot type containing a positive-measure set of invariant tori. However, when considered as magnetohydrostatic equilibria, Beltrami fields have constant magnetic pressure, so this construction does not provide a foliation by regular pressure surfaces. For nonconstant pressure, Bruno and Laurence~\cite{BL} proved existence of piecewise smooth equilibria with pressure jumps in toroidal domains close to axisymmetric ones. A subsequent construction~\cite{ELP} gives piecewise smooth equilibria with arbitrarily many current sheets in toroidal domains of any knot type. Their pressure is piecewise constant, and the field and pressure need not be continuous across the interfaces. In \cite{EPP} we constructed counterexamples to Grad's conjecture in toroidal domains of arbitrary knot type using a new convex integration scheme, but the magnetic field is only H\"older continuous.

In order to state our main theorem, let us introduce some notation. A magnetohydrostatic equilibrium in a bounded domain $\Omega\subset\R^3$ is a magnetic field $B$ and a pressure function $P$ satisfying
\begin{equation}
\Div B = 0, \qquad (\curl B) \times B = \nabla P \quad \text{in }\Omega, \qquad B \cdot n = 0 \quad \text{on }\partial\Omega.
\label{kg2}
\end{equation}
Here the differential operators and the cross product are Euclidean, and $n$ is the outward unit normal. Taking the scalar product of the force equation with $B$ shows that the pressure is constant along magnetic lines (in other words, $P$ is a first integral of $B$).

As required in Grad's conjecture, we are interested in equilibria with nested pressure tori surrounding a single nondegenerate magnetic axis. Our goal is to prescribe the knot type of this entire configuration while keeping the field nonzero throughout the domain, including its axis and boundary. The equilibria we aim to construct have no nontrivial Euclidean symmetry (discrete or continuous) and occur in smooth families of pairwise inequivalent configurations. Their magnetic flows can have any prescribed positive slope (i.e., constant rotational transform) in suitable angular coordinates.

Our starting point is the smooth axisymmetric equilibria constructed in~\cite{PSS}. They realize Seifert fibrations with prescribed coprime nonzero winding numbers in rotationally symmetric toroidal domains, and their construction also allows irrational trajectories. The ambient solid torus and its magnetic axis remain unknotted. Our construction deforms these equilibria to prescribe the knot type of the axis and the entire pressure foliation. As in~\cite{GSLT}, the domain boundary is part of the construction. The new geometric idea that we develop in this work is to deform a finite covering immersion of the aforementioned axisymmetric equilibria, separating its coincident sheets into a prescribed braid while preserving force balance and magnetic dynamics. This allows changes of knot type via Alexander's theorem in knot theory, which cannot be achieved by deforming an embedding. We believe that this strategy may also be useful for other differential equations that admit a deformation theory on immersed domains.

The following is our main result. We recall that the knot type of a solid torus is that of its core. We consider smooth knot types in $\R^3$ and write $\T = \R/(2\pi\Z)$.

\begin{theorem}\label{kg3}
For every knot type and real~$\kappa > 0$, there is a magnetohydrostatic equilibrium $(B, P)$ solving~\eqref{kg2} on a smoothly embedded solid torus~$\Omega$ of that knot type. This equilibrium is of class~$C^\infty$ up to the boundary and $B$ does not vanish at any point of $\overline\Omega$. Furthermore:
\begin{enumerate}
\item The pressure is positive in $\Omega$ and zero on $\partial\Omega$. Its critical set is a single embedded core curve of $\Omega$, with negative definite normal Hessian, and all its regular nonempty level sets are embedded tori.
\item On each regular pressure torus, $B$ is smoothly conjugate, with time unchanged, to
\begin{equation}
\omega\left(\partial_\Theta + \kappa\partial_\zeta\right) \quad \hbox{on }\T^2,
\label{kg4}
\end{equation}
where $\omega \equiv \omega(P) > 0$ depends smoothly on the pressure level. Here, $(\Theta, \zeta)$ are suitable angular coordinates parametrizing $\mathbb T^2$.
\end{enumerate}
\end{theorem}

Theorem~\ref{kg3} immediately implies that the result also holds for all $\kappa<0$.
Indeed, let us fix a Euclidean reflection $\mathcal R$. We apply the theorem with slope~$-\kappa > 0$ to the mirror of the prescribed knot type and denote by $(B, P)$ the resulting equilibrium on $\Omega$. Since the orientation signs in curl and cross product cancel, we see that $\widetilde B = \mathcal R_* B$ and $\widetilde P = P \circ \mathcal R^{-1}$ solve~\eqref{kg2} on $\widetilde\Omega = \mathcal R(\Omega)$, which has the prescribed knot type. On each regular pressure torus, we use the angular coordinates~$\widetilde\Theta = \Theta \circ \mathcal R^{-1}$ and $\widetilde\zeta = -\zeta \circ \mathcal R^{-1}$. Then, we have
\begin{equation*}
\widetilde B = \omega(\widetilde P) \left(\partial_{\widetilde\Theta} + \kappa\partial_{\widetilde\zeta}\right).
\end{equation*}
The longitudinal reversal compensates for the ambient orientation reversal. The conjugacy preserves time, and all other conclusions are preserved by the reflection.

For rational $\kappa$, all magnetic lines are periodic. For irrational $\kappa$, every magnetic line is dense in its pressure torus and the flow there is uniquely ergodic (including for Liouville slopes). The construction provides a smooth radial label $\varrho \geq 0$ on $\overline\Omega$, vanishing precisely on the axis, such that $P = p(\varrho)$ with $p' < 0$. By composition, we can write the speed factor in~\eqref{kg4} as $\omega(\varrho)$.

The knot representative and the domain boundary are selected by the construction. The angles in~\eqref{kg4} give a meridian and a longitude transported from the covered source solid torus. The latter can differ from the preferred zero-linking longitude by an integral meridian; Section~\ref{kg106} gives the corresponding change of slope. To obtain any prescribed $\kappa>0$, we choose a covering degree~$N$ so that the reference axisymmetric equilibrium has slope $\kappa_* \coloneqq N\kappa>1$. The value of $\kappa_*$ plays an important role in the proof of coercivity, while the flow on the cover has slope $\kappa$.

\begin{remark}\label{kg5}
A variation of this construction gives equilibria of any prescribed finite regularity~$C^r$ with nonconstant winding (that is, rotational transform), so that rational and irrational pressure tori coexist; see Proposition~\ref{kg116}. The admissible change in rotational transform depends on the regularity through highest-order commutator estimates. As a result, our argument does not yield a $C^\infty$ equilibrium in this case.
\end{remark}

The construction also gives geometric freedom within each knot type. Two equilibria are equivalent under similarities and field rescalings if a Euclidean similarity maps one domain onto the other and pushes one field to a nonzero constant multiple of the other. The force equation then determines the pressure transformation up to an additive constant.

\begin{corollary}\label{kg6}
The equilibria in Theorem~\ref{kg3} can be chosen in a smooth one-parameter family whose members have no nontrivial Euclidean symmetries, even allowing reversal of the magnetic field. Distinct members are inequivalent under similarities and constant field rescalings.
\end{corollary}

Smoothness of a family on varying domains is understood after pullback by a smooth family of embeddings of a fixed solid torus. Such a common source domain is available locally around every member of the family in the corollary. Note that, for a nontrivial knot, the absence of continuous Euclidean symmetry already follows from its critical axis. Indeed, a closed embedded curve that is invariant under rotations around some axis must itself be a circle.

The article is organized as follows. In Section~\ref{kg7} we present the main ideas of the proof of the main theorem. Section~\ref{sec reference} constructs the reference equilibrium, and Section~\ref{kg29} derives the nonlinear equation. In Section~\ref{kg41} we prove the linear inverse estimates and Section~\ref{kg90} constructs and integrates the braiding displacement. Finally, in Section~\ref{kg106} we complete the geometric argument. The case of nonconstant rotational transform at finite regularity, cf. Remark~\ref{kg5}, is considered in Section~\ref{kg115}.

\section{Ideas of the proof}\label{kg7}

The main unknown in our construction is an immersion $\Phi$ of a fixed
solid torus into $\R^3$, which we use to transport a fixed magnetic field from the solid torus. To construct magnetic equilibria, our goal is of course to find immersions that preserve force balance (i.e., the equilibrium equation), which leads to a number of technical challenges that we will next sketch. Before getting into technicalities, however, let us emphasize that the simplest yet perhaps the most potentially useful idea of the proof is that one can attempt  (and in fact manage) to implement this strategy when $\Phi$ is an immersion, instead of an embedding as is done in the deformation arguments that we have seen in the literature. In our proof, it is essential to choose the immersion so that the magnetic axis becomes a knot of the prescribed type. 

Let us now pass to a sketch of how the proof goes on. For this, we first state the nonlinear equation for $\Phi$ and the linear equation used to solve it.

\medskip\noindent\emph{The source torus and its magnetic field.} By Alexander's theorem~\cite{Alexander}, we may represent the prescribed knot as a closed braid. Stabilizing this braid increases the number of strands $N \geq 2$, so we can arrange that $\kappa_* := N\kappa > 1$. We then use the construction in~\cite{PSS} to obtain an axisymmetric reference equilibrium of slope~$\kappa_*$, depending smoothly on its thickness~$\sigma > 0$ in a fixed interval. Its lift to the $N$-fold longitudinal cover has covering immersion~$\Phi_0$. We will deform this immersion to separate the coincident sheets according to the chosen braid. Specifically, the key geometric observation is that the $N$ sheets
of the covering immersion coincide initially, but
can move independently under a deformation of the
source. At a fixed longitudinal angle $\theta$, their
first-order transverse positions are prescribed by
the $N$ distinct points
\[
\mathbf F(\theta+2\pi j), \qquad j=0,\ldots,N-1.
\]
Thus, if the deformation realizes this first variation,
the coincident sheets separate into a closed braid.
Although the full immersion need not be injective,
a sufficiently small neighborhood of the resulting
embedded axis will be.

The domain of the immersion, which we will call the {\em source solid torus}, is the fixed manifold
\begin{equation*}
\Omega_{\sigma, N} \coloneq \{(R, z): R > 0, \;(\log R)^2 + z^2 < \sigma^2\} \times \T_N, \qquad \T_N = \R/(2\pi N\Z).
\end{equation*}
The initial map is the covering immersion $\Phi_0(R, z, \varphi) \coloneqq (R\cos\varphi, R\sin\varphi, z)$. We seek a deformation $\Phi: \overline{\Omega}_{\sigma, N} \to \R^3$ of this map, keeping the magnetic field~$b$ and the volume form~$\mu$ fixed on $\Omega_{\sigma, N}$. Thus the domain of the map stays fixed, while its image varies with $\Phi$. To obtain the physical domain, we will restrict $\Phi$ to a smaller solid torus on which it is an embedding.

On this source torus we use a squared radial coordinate~$\varrho$, which vanishes on the core axis, and angular coordinates~$\Theta$ of period~$2\pi$ and $Z$ of period~$2\pi N$. In these coordinates, the reference magnetic field, its current and the volume form take the form
\begin{equation}
b = \omega(\varrho)(\partial_\Theta - \kappa_* \partial_Z), \qquad j_0 = \beta_0(\varrho) \partial_Z, \qquad \mu = d\varrho \wedge d\Theta \wedge dZ,
\label{kg8}
\end{equation}
with $\omega, \beta_0 > 0$. As just explained, both the field~$b$ and the form~$\mu$ are kept fixed during the deformation. The change of coordinates preserves physical time, and using the angle $\zeta = -Z/N$ gives the slope~$\kappa_*/N = \kappa$ required by the theorem. The reference equilibrium and the adapted coordinates are constructed in Section~\ref{sec reference}; see in particular Lemma~\ref{kg13} and~\eqref{kg20}--\eqref{kg23}.

\medskip\noindent\emph{The nonlinear equation for $\Phi$}. To write the equation on the fixed domain, for a one-form~$\eta$ on $\Omega_{\sigma, N}$ we define the field~$\curl_\mu \eta$ by $\iota_{\curl_\mu \eta} \mu = d\eta$. Given a candidate immersion~$\Phi$, let
\begin{equation*}
G_\Phi := \Phi^* g_{\rm Eucl}, \qquad \Phi^* dV =: \rho_\Phi \mu, \qquad c_\Phi := \curl_\mu\bigl(G_\Phi(b, \cdot)\bigr),
\end{equation*}
where $dV$ is the Euclidean volume form. The quantities $G_\Phi$, $\rho_\Phi$, and $c_\Phi$ are all computed from $\Phi$, so they introduce no additional unknowns. If $\rho_\Phi = 1$, then $c_\Phi$ is the Euclidean current pulled back to $\Omega_{\sigma, N}$. For an arbitrary candidate immersion it is still divergence-free: this follows from defining curl with the fixed form~$\mu$ on $\Omega_{\sigma, N}$.

We denote by $\av$ the average of vector fields on $\Omega_{\sigma, N}$ under poloidal rotations in $\Theta$ and longitudinal translations in $Z$. The nonlinear equation we impose is
\begin{equation}
\mathcal F(\Phi) \coloneqq \left(\rho_\Phi - 1, (\Id - \av) c_\Phi\right) = 0.
\label{kg9}
\end{equation}
The first component requires the immersion to preserve volume; the second requires $c_\Phi$ to be invariant under the two angular actions. Both conditions are satisfied by the reference covering immersion~$\Phi_0$. Volume preservation ensures that the transported magnetic field is divergence-free. To understand the second condition, note that, by Lemma~\ref{lem:profiles}, a divergence-free current invariant under both angular actions must have the form
\begin{equation*}
c_\Phi = \alpha(\varrho) b + \beta(\varrho) \partial_Z
\end{equation*}
where the profiles $\alpha, \beta$ are obtained from $c_\Phi$. If the immersion is a sufficiently small deformation of the reference map, these profiles remain close to those of the reference current, and in particular $\beta > 0$. It follows that the transported magnetic field solves~\eqref{kg2} on $\Phi(\Omega_{\sigma, N}) \subset \R^3$, with field and pressure given by
\begin{equation*}
B \circ \Phi = D\Phi\,b, \qquad P \circ \Phi = p(\varrho), \qquad p' = -\beta\omega.
\end{equation*}
Indeed, the one-form associated with $(\curl B) \times B$ pulls back to $-\beta\omega\,d\varrho$, which is precisely the pressure differential above. Adding a constant makes $P$ vanish on the boundary. We impose no boundary condition on the displacement, so this boundary is free to move. Solving $\mathcal{F}(\Phi) = 0$ therefore produces an MHD equilibrium; the remaining geometric task is to choose a solution $\Phi$ that separates the $N$-fold covered magnetic axis into the prescribed braid. The precise correspondence between solutions of~\eqref{kg9} and magnetohydrostatic equilibria is established in Proposition~\ref{kg33}.

\medskip\noindent\emph{The linearized problem.} Fix an immersion~$\Phi$ close to $\Phi_0$ and set~$j \coloneqq \av c_\Phi$. This averaged current has the radial form~$j = \alpha(\varrho) b + \beta(\varrho) \partial_Z$. The linear system that we shall solve is
\begin{align}
\Div_\mu \xi &= h,
\label{kg10}
\\
(\Id - \av) \left\{\curl_\mu\left(G_\Phi([b, \xi], \cdot)\right) - [j, \xi] + (\Div_\mu \xi) j\right\} &= g.
\label{kg11}
\end{align}
The unknown in this system is the displacement field~$\xi$ on $\Omega_{\sigma, N}$; the coefficients~$b, G_\Phi, j$ are fixed. We prescribe a scalar datum~$h$ and a vector datum~$g$, with $\Div_\mu g = 0$ and $\av g = 0$. Let $M_\Phi$ denote the operator on the left-hand side.

At the reference immersion~$\Phi_0$, this system is exactly the linearization of~\eqref{kg9}. To see this, represent an infinitesimal variation of the immersion as $\dot\Phi = D\Phi_0\,\xi$, where $\xi$ is a vector field on the source. Direct differentiation gives
\begin{equation*}
D\mathcal F(\Phi_0)[D\Phi_0\,\xi] = M_{\Phi_0} \xi.
\end{equation*}
For a nearby immersion~$\Phi$, the expressions $D\mathcal F(\Phi)[D\Phi\,\xi]$ and $M_\Phi \xi$ differ by an error bilinear in $\mathcal F(\Phi)$ and $\xi$. This is the error allowed by the Nash--Moser theorem used below. The exact linearization formula and the bilinear error identity are given in Proposition~\ref{kg36} and~\eqref{kg37}.

Our analytic task is consequently to construct a smooth tame family of right inverses: a map~$R_\Phi$ which, for each datum~$(h, g)$, gives a solution $\xi = R_\Phi(h, g)$ of~\eqref{kg10}-\eqref{kg11} with suitable dependence on the data. We first subtract a particular solution of the divergence
equation, leaving a divergence-free displacement, and absorb
its contribution to~\eqref{kg11} into the datum~$g$.
Introducing $w \coloneqq [b,\xi]$ and imposing $\av\xi=0$,
we are left with the coupled system
\begin{align*}
w &= [b,\xi], \qquad \Div_\mu\xi=0,\\
[j,\xi]
&= (\Id-\av)\curl_\mu\bigl(G_\Phi(w,\cdot)\bigr)-g.
\end{align*}
Here $w$ is an unknown, not a
prescribed magnetic variation. We do not need to invert
magnetic transport on arbitrary vector fields: we must
solve the coupled system, which will select a $w$ compatible
with that transport.

\medskip\noindent\emph{Solving the linear system.} The equation $w = [b, \xi]$ differentiates the displacement along magnetic lines. Recovering $\xi$ directly is difficult at Fourier modes which vary very little along $b$. Indeed, for the field in parentheses in~\eqref{kg8} and a scalar angular mode~$(n, m) \in \Z \times N^{-1} \Z$, we have
\begin{equation*}
(\partial_\Theta - \kappa_* \partial_Z) e^{\ii(n\Theta+mZ)} = \ii(n - \kappa_* m) e^{\ii(n\Theta+mZ)} \equiv \ii k\,e^{\ii(n\Theta+mZ)}.
\end{equation*}
Thus $k = 0$ means that the scalar mode is constant along the magnetic lines. Even when $k \neq 0$, inversion requires division by $k$, which can lead to small divisors. We therefore perform this inversion only for $\abs{k} \geq \delta$, where $\delta > 0$ is a fixed small number. The corresponding spectral decomposition and the off-band inverse are constructed in~\eqref{kg46} and Lemma~\ref{kg47}.

For the nonzero angular modes in the complementary band $\abs{k}<\delta$, we retain the radial displacement
\begin{equation*}
f \coloneqq \omega\,d\varrho(\xi)
\end{equation*}
as an additional unknown, with $d\varrho(w)=\ii kf$. Crucially, $w$ is itself an unknown, to be determined together with $f$, rather than an arbitrary prescribed magnetic variation. Keeping $f$ avoids division by $k$, but the full relation $w=[b,\xi]$ must still be enforced. At exact resonance, for instance, this relation requires
\[w=-\frac{\omega'}{\omega}f(\partial_\Theta-\kappa_*\partial_Z), \qquad f=\frac{\omega}{\kappa_*\omega'}\,w^Z.\]
Thus $f$ remains algebraically constrained even when $d\varrho(w)=0$, since $\omega'=6+O(\sigma)>0$; see~\eqref{kg23}. Here primes denote differentiation with respect to $\varrho$.

We impose $\av\xi=0$, hence the corresponding zero-average
condition on $(w,f)$, and require $w$ to be divergence-free.
Once $w$ is known, we reconstruct $\xi$ using the current
equation
\begin{equation*}
[j,\xi]
=(\Id-\av)\curl_\mu\bigl(G_\Phi(w,\cdot)\bigr)-g.
\end{equation*}
The frequency of current transport is
$\lambda=\alpha\omega k+\beta m$, which equals
$\beta_0m$ at the reference equilibrium.
Choosing $\delta>0$ sufficiently small ensures that
all selected modes have $n,m\neq0$ and, for nearby
currents,
\begin{equation*}
\abs{\lambda}\geq c\abs{m}\geq c/N.
\end{equation*}
Thus current transport allows reconstruction without
division by $k$, even at exact magnetic resonance.
The variational equations described next determine
$(w,f)$ and ensure that the reconstructed displacement
also satisfies $w=[b,\xi]$, including the algebraic
compatibility at $k=0$; see Lemma~\ref{kg60}.

To find $w$ and $f$ from the datum~$g$, we put the current
equation in Hermitian variational form on
$M=\Sigma\times\T_N$, extending the transverse disk to a
sphere $\Sigma$ to avoid imposing boundary conditions.
The zero-angular-average condition removes the harmonic
fields and allows us to invert curl.

The key is the sign structure of the form~\eqref{kg59}.
Writing $\tau$ for the selected-band component of $w^Z$,
its scalar terms satisfy
\[\langle\mathsf Cf,f\rangle +2\operatorname{Re}\langle\mathsf B\tau,f\rangle = \bigl\|\mathsf C^{1/2} (f+\mathsf C^{-1}\mathsf B\tau)\bigr\|_{L^2}^2-\bigl\langle\mathsf B^2\mathsf C^{-1}\tau,\tau \bigr\rangle,\]
where $\mathsf B,\mathsf C$ are the bounded real modal
coefficients of~\eqref{kg56}, with $\mathsf C\geq c>0$.
For $k\neq0$, eliminating $f$ would insert
$w^\varrho/(\ii k)$ inside the positive square.
Thus the potentially singular contribution has a
favorable sign and can be kept on the coercive side
of the estimate. Retaining $f$ makes this formulation
valid also at exact resonance.

The nontrivial step is to absorb the negative term
into the metric energy. Its coefficient is close to
$1/(3\kappa_*^2)$, while the solenoidal estimate of
Lemma~\ref{kg67} bounds $\|\tau\|_{L^2}^2$ by the full
metric energy with leading constant $1+2\kappa_*^2$.
This estimate uses the divergence-free constraint
and retains the shear of the reference metric;
a pointwise metric comparison would not suffice.
The coefficient left in front of the metric energy
is therefore
\begin{equation*}
1-\frac{1+2\kappa_*^2}{3\kappa_*^2}
=\frac{\kappa_*^2-1}{3\kappa_*^2}>0,
\end{equation*}
up to controlled small errors. This margin allows us
to absorb the remaining perturbative terms.
Together with the positive square, the metric energy
then controls both $w$ and $f$, without a lower bound
for the nonzero values of $\abs{k}$ inside the band.
Proposition~\ref{kg69} establishes this coercivity,
treating exact resonances and Liouville slopes
within the same argument.

The Lax--Milgram theorem yields a weak solution of the
augmented variational problem. To prove smoothness, we commute the operator with all powers of the appropriate Sobolev operator. The commutators have lower order in the unknowns, so their contribution is controlled by interpolation. In this way, estimates at every order hold for one fixed thickness and coefficient neighborhood, chosen using only finitely many low-order norms. Combining these estimates with the reconstruction of $\xi$ gives the inverse $R_\Phi$, with smooth dependence and a fixed loss of derivatives. The higher-order estimates and smooth tame dependence are proved in Proposition~\ref{kg75}, and the right inverse is assembled in the proof of Proposition~\ref{prop right inverse of the linear operator}.

\medskip\noindent\emph{Choosing a linear solution that creates the knot.} We shall use this inverse to find a solution of the homogeneous linearized equation whose displacement moves the core axis along the prescribed braid. First, we subtract the common mean of the strands in each section, thereby centering the braid without changing strand separation. We then approximate the centered profile by a trigonometric polynomial. Centering removes the integer longitudinal frequencies, and the approximating profile contains only finitely many frequencies~$m$.

We also choose the profile so that $\kappa_* m \neq \pm1$ at each of its frequencies. The reason is that, on the axis, the map from normal displacement to magnetic variation has determinant proportional to $1 - \kappa_*^2 m^2$. Its invertibility will turn a small magnetic correction into a small correction of the axis. Since $m = \ell/N$ with $\ell \in \Z$, the frequencies to avoid are those satisfying $\kappa\ell = \pm1$. If such modes occur, we remove them by applying a common periodic rotation of winding number zero to the normal sections and then truncating the Fourier series. These operations preserve the knot type. For the existence of a centered trigonometric braid profile satisfying these nonresonance conditions, see Lemma~\ref{kg93}.

At the reference immersion the inverse preserves each longitudinal frequency. For the finite set of frequencies just chosen, we take $\delta$ so that the selected band contains only the modes with $k = 0$. The differentiated radial terms of the augmented variational system then vanish. Moreover, its forcing is represented by a vector potential obtained through radial integration, which contributes a factor $\sigma$. The uniform estimates after rescaling the cross-section to unit size therefore yield an $O(\sigma)$ correction of the magnetic field. By the invertibility of the axis map, the normal correction on the axis is also $O(\sigma)$, as shown in Corollary~\ref{kg88}.

We extend the braid profile to a solenoidal trial displacement $\xi_{\rm tr}$, whose defect $M_{\Phi_0} \xi_{\rm tr}$ in the linear equation stays bounded in these rescaled norms. We can now correct it by setting
\begin{equation*}
\xi_{\rm corr} = R_{\Phi_0}(M_{\Phi_0} \xi_{\rm tr}), \qquad \xi = \xi_{\rm tr} - \xi_{\rm corr}.
\end{equation*}
This gives $M_{\Phi_0} \xi = 0$, while the normal motion defined by $\xi$ differs from the trial braid by $O(\sigma)$. The order of choices is therefore to fix the braid, then the band, and finally take $\sigma$ small enough that strand separation and the knot type are preserved. The details are worked out in~\eqref{kg99}--\eqref{kg101}.

\medskip\noindent\emph{From the linear solution to exact equilibria.} We now look for exact solutions whose initial displacement is $D\Phi_0\,\xi$. The controlled error between $M_\Phi$ and the true derivative permits the use of Hamilton's implicit-function theorem with quadratic error~\cite[Part III, Theorem 3.3.1]{Hamilton}. To ensure that the initial displacement is the one we have chosen, write
\begin{equation*}
\Phi = \Phi_0 + t(D\Phi_0\,\xi + v),
\end{equation*}
where $t$ is the deformation parameter and $v$ is the unknown correction to the immersion. We apply Hamilton's theorem to the equation $\mathcal F(\Phi)/t = 0$, which extends smoothly to $t = 0$. At that parameter value it is solved by $v = 0$, since $M_{\Phi_0} \xi = 0$. The solution supplied by the theorem satisfies $v(t) = O(t)$, and hence
\begin{equation*}
\mathcal F(\Phi_t) = 0, \qquad \Phi_t = \Phi_0 + tD\Phi_0\,\xi + O(t^2).
\end{equation*}
It follows that the normalized transverse motion of the axis tends to the corrected braid profile. For sufficiently small~$t > 0$, the axis is therefore an embedded knot of the prescribed type. We restrict $\Phi_t$ to a sufficiently small tube in the source, obtaining an embedded solid torus. The formulas for $B$ and $P$ then give the equilibrium, its pressure foliation and its magnetic dynamics, with physical time unchanged. The passage from these exact immersions to embedded equilibria, including the prescribed magnetic dynamics, is completed in Section~\ref{kg106}.

To prove Corollary~\ref{kg6}, we make some additional choices in the braid profile. Any symmetry of an equilibrium preserves the unique critical pressure axis, even if it reverses the field. We arrange for this axis to have a unique maximum of curvature at which the third arclength derivative of curvature is nonzero. A symmetry must fix that point; it must also preserve the direction along the axis, and hence fix the axis pointwise. These conditions are arranged in Lemma~\ref{kg95}, and the resulting rigidity is established in Lemma~\ref{kg110}. Positive curvature and nonzero torsion at the same point then imply that the ambient isometry is the identity. To distinguish the members of the family, we use the product of the axis length and the integral of the squared deviation of curvature from its arclength average. This quantity is invariant under similarities and varies strictly along the family; see~\eqref{kg113}--\eqref{kg114}. Since an equivalence of equilibria maps their pressure axes to one another, distinct members are inequivalent under similarities and constant field rescalings.

\section{The reference equilibrium and the adapted coordinates}
\label{sec reference}

Our purpose in this section is to construct the axisymmetric equilibrium from which the deformation will start. Its magnetic lines have constant slope on the pressure tori. After lifting the equilibrium to a finite cover, we shall introduce coordinates in which the volume form is standard and the coefficients of both the magnetic field and its current depend only on the radial variable. This covered solid torus will be the domain of our nonlinear equation. We shall also extend its reference coefficients to a closed manifold, as needed to construct the right inverse.

The starting point is the smooth planar construction of~\cite[Proposition~2.1]{PSS}. We state below the corresponding family on the disk, together with the first-order expansion that enters our estimates. We stress that the parameter interval is the same at all differentiability orders.

\begin{lemma}\label{kg13}
Let $D \subset \R^2$ be the unit disk and $\Psi_0(x_1, x_2) = (x_1^2 + x_2^2)/2$. There are a constant $\sigma_0 > 0$, and area-preserving diffeomorphisms $\gamma_\sigma: \overline D \to \overline D$ and functions~$F_\sigma \in C^\infty([0, 1/2])$ for $\abs{\sigma} < \sigma_0$, which satisfy the following (in the image coordinates~$(x, y) = \gamma_\sigma(x_1, x_2)$):
\begin{equation}
\Psi_\sigma = \Psi_0 \circ \gamma_\sigma^{-1}, \qquad L_\sigma \Psi_\sigma = F_\sigma(\Psi_\sigma), \qquad L_\sigma = e^{-4\sigma x}(\partial_{xx} - 2\sigma\partial_x) + e^{-2\sigma x} \partial_{yy}.
\label{kg14}
\end{equation}
Both families are jointly smooth in their variables and in $\sigma$, up to the boundary. Moreover,
\begin{equation}
\gamma_\sigma(x_1, x_2) = (x_1, x_2) + \sigma(1 - x_1^2 - 3x_2^2, 2x_1 x_2) + O(\sigma^2), \qquad F_\sigma = 2 + O(\sigma^2),
\label{kg15}
\end{equation}
where the remainders divided by $\sigma^2$ extend jointly smoothly to $\sigma = 0$, with uniform bounds in each fixed smooth seminorm on compact parameter subintervals.
\end{lemma}

\begin{proof}
The existence of the area-preserving family, as well as its joint smoothness, is established in~\cite[Proposition~2.1 and its proof]{PSS}. In that argument the parameter is allowed to vary in an interval $(-\sigma_0,\sigma_0)$, so the family is also defined at $\sigma = 0$. We choose the normalization so that the planar curl of $\gamma_\sigma - \Id$ has zero angular average. With this choice, the first-variation computation in~\cite[Section~3]{PSS} yields the expansion of $\gamma_\sigma$ stated in~\eqref{kg15}. To check the joint smoothness of $F_\sigma$, notice that $(L_\sigma \Psi_\sigma) \circ \gamma_\sigma$ is jointly smooth and radial on the disk. It is therefore a smooth function of $(\sigma, (x_1^2 + x_2^2)/2)$ even at the center, as required.

For completeness, let us compute the first derivative of $F_\sigma$ with respect to $\sigma$ directly from its equation. In the image coordinates, the expansion of $\gamma_\sigma$ implies that
\begin{equation*}
\left.\partial_\sigma \Psi_\sigma\right|_0 = x(x^2 + y^2 - 1), \qquad \Delta\bigl(x(x^2 + y^2 - 1)\bigr) = 8x, \qquad \left.(\partial_\sigma L_\sigma)\right|_0 \Psi_0 = -8x.
\end{equation*}
Differentiating~\eqref{kg14} with respect to the parameter and using $F_0 = 2$, we obtain $\left.\partial_\sigma F_\sigma\right|_0 = 0$. The asserted bounds for the remainders then follow by applying Taylor's formula to the jointly smooth families, in any fixed smooth seminorm and on the same parameter interval.
\end{proof}

We now turn the planar solution into an axisymmetric equilibrium on a finite cover. For this, we fix a parameter \(\kappa_*>1\), which determines the slope of the reference equilibrium. In the proof of Theorem \ref{kg3}, we will choose the covering degree \(N\), together with the braid, so that \(\kappa_*=N\kappa\), where \(\kappa>0\) is the prescribed slope in the theorem. For the moment, given an integer~$N \geq 2$, we denote $\T_N = \R/(2\pi N\Z)$ and $\T = \T_1$. The covered solid torus on which we work is
\begin{equation}
\Omega_{\sigma, N} = \{(R, z): (\log R)^2 + z^2 < \sigma^2\} \times \T_N
\label{kg16}
\end{equation}
and the reference immersion is defined by
\begin{equation}
\Phi_0(R, z, \varphi) = (R\cos\varphi, R\sin\varphi, z), \qquad G_0 = \Phi_0^* g_{\rm Eucl}, \quad \mu = \Phi_0^* dV.
\label{kg17}
\end{equation}
This map has rank three at every point and covers its image $N$ times. We use the lifted cylindrical frame to define the fields and differential operators on the source domain. Throughout this construction, unless stated otherwise, the curl and gradient are those of $G_0$, while the divergence is computed with respect to $\mu$.

Using the rescaled meridional coordinates~$(x, y) = (\log R/\sigma, z/\sigma)$, we define the stream function and the magnetic field as
\begin{equation}
\psi = \sigma^2 \Psi_\sigma(x, y), \qquad b = -\frac{\partial_z\psi}{R} e_R + \frac{\kappa_*}{R} e_\varphi + \frac{\partial_R\psi}{R} e_z.
\label{kg18}
\end{equation}
This vector field is divergence-free and does not vanish anywhere. To see that it is an equilibrium, observe that the above rescaling transforms the Grad--Shafranov equation into~\eqref{kg14}. Moreover, in the coordinates~$(x, y) = \gamma_\sigma(x_1, x_2)$ we have $\psi = \sigma^2 \mathcal I$, with $\mathcal I = (x_1^2 + x_2^2)/2$. Its current is toroidal, and the force is given by $(\curl b) \times b = -F_\sigma(\mathcal I) \nabla\psi$. We therefore obtain an exact equilibrium by taking the pressure to be $\sigma^2 \int_{\mathcal I}^{1/2} F_\sigma(s)\,ds$.

It remains to straighten the flow without changing its physical time. Let us denote by $J(x_1, x_2) = (-x_2, x_1)$ the rotation through a right angle. Since the planar change of coordinates preserves area, it satisfies $(\gamma_\sigma)_* \partial_\theta = J\nabla\Psi_\sigma$, where $\theta = \arg(x_1 + \ii x_2)$. Hence, setting $Z_0 = -\varphi$, the field and the volume form read as
\begin{equation*}
b = R^{-2}(\partial_\theta - \kappa_* \partial_{Z_0}), \qquad \mu = \sigma^2 R^2\,d\mathcal I \wedge d\theta \wedge dZ_0.
\end{equation*}
In these formulas, $\mathcal I$ denotes the enclosed area in the rescaled disk divided by $2\pi$; the cylindrical radius~$R$ is now regarded as a function of $(\mathcal I, \theta)$. We introduce the new variables through the following equations:
\begin{equation}
\begin{gathered}
\omega(\mathcal I) = \langle R^2\rangle_\theta^{-1}, \qquad \partial_\theta u = \omega(\mathcal I) R^2 - 1, \qquad \langle u\rangle_\theta = 0, \\
\Theta = \theta + u, \qquad Z = Z_0 - \kappa_* u, \qquad \varrho(\mathcal I) = \sigma^2 \int_0^{\mathcal I} \omega(s)^{-1} \dd s.
\end{gathered}
\label{kg19}
\end{equation}
Here and in what follows, \(\langle\cdot\rangle_\theta\) denotes the average over \(\theta\) at fixed \(\mathcal I\).

Let us check the regularity of the angular primitive at the centre. The angular primitive with zero mean is obtained by integration against a bounded kernel for the rotation action, an operation which preserves smoothness in $(x_1, x_2)$. The angular change is also nondegenerate, since $1 + u_\theta = \omega R^2 > 0$. For the radial variable, the formula $\varrho/\mathcal I = \sigma^2 \int_0^1 \omega(t\mathcal I)^{-1}\,dt$ shows that this quotient has a smooth positive extension to zero. We may therefore define a Cartesian coordinate~$y = (y_1, y_2)$ on the disk by
\begin{equation*}
y = \sqrt{\varrho/\mathcal I}\,R_{u(x_1, x_2)}(x_1, x_2), \qquad \varrho = \tfrac12 \abs{y}^2,
\end{equation*}
where $R_u$ denotes planar rotation. Together with the longitudinal coordinate~$Z = Z_0 - \kappa_* u$, this gives smooth coordinates on the whole closed source torus, including the axis. Notice that the term added to $Z_0$ is a single-valued real function. The new longitudinal coordinate is thus well defined modulo $2\pi N$ for every real $\kappa_*$, with no arithmetical restriction on the slope.

Since the change $\mathcal I \mapsto \varrho$ is increasing, it does not alter the level sets. Its role is to put the volume form in product form with respect to the straightened angles. We shall henceforth regard $\omega$ as a function of $\varrho$ and write $X = \partial_\Theta = Jy \cdot \partial_y$ and 
$$U = \partial_Z\,.$$ 
With this notation, the preceding construction gives
\begin{equation}
V = X - \kappa_* U, \qquad b = \omega(\varrho) V, \qquad \mu = d\varrho \wedge d\Theta \wedge dZ = dy_1 \wedge dy_2 \wedge dZ,
\label{kg20}
\end{equation}
on $0 \leq \varrho \leq \varrho_b$, where
\begin{equation}
\varrho_b = \sigma^2 \int_0^{1/2} \omega(\mathcal I)^{-1} \dd \mathcal I = \tfrac12\sigma^2(1 + O(\sigma^2)).
\label{kg21}
\end{equation}
We shall use $Y^\varrho = d\varrho(Y)$ to denote the radial component of a vector field~$Y$. With the coordinates just constructed, the current and the force relation can be written as
\begin{equation}
j_0 = \curl b = F_\sigma(\mathcal I) U =: \beta_0(\varrho) U, \qquad b \times U = \nabla\psi = \omega\nabla \varrho.
\label{kg22}
\end{equation}
In particular, $b$ and $j_0$ preserve the volume form~$\mu$ and commute with one another.

We shall need the leading terms of these coefficients to control the sign of the quadratic form in the linear problem. From the expansion~\eqref{kg15}, one obtains
\begin{equation*}
R^2 = 1 + 2\sigma x_1 + \sigma^2(2 - 6x_2^2) + O(\sigma^3), \qquad \omega(\mathcal I) = 1 - 2\sigma^2 + 6\sigma^2 \mathcal I + O(\sigma^3).
\end{equation*}
To differentiate the radial profiles, we use~\eqref{kg19}, which gives $\partial_\varrho = (\omega/\sigma^2) \partial_{\mathcal I}$. Accordingly, if primes denote derivatives with respect to $\varrho$, the estimates take the form
\begin{equation}
\omega = 1 + O(\sigma^2), \qquad \omega' = 6 + O(\sigma), \qquad \beta_0 = 2 + O(\sigma^2), \qquad \abs{\beta_0'} \leq C.
\label{kg23}
\end{equation}
The bounds for the radial derivatives deserve a brief justification. Let us write
\begin{equation*}
F_\sigma(\mathcal I) = 2 + \sigma^2 r_F(\sigma, \mathcal I), \qquad \omega(\mathcal I) = 1 + \sigma^2(6\mathcal I - 2) + \sigma^3 r_\omega(\sigma, \mathcal I),
\end{equation*}
where the two remainders are jointly smooth for the parameters in a common interval and $0 \leq \mathcal I \leq 1/2$. Differentiating these expressions yields
\begin{equation*}
\beta_0' = \omega\,\partial_{\mathcal I} r_F, \qquad \omega' = \omega\bigl(6 + \sigma\partial_{\mathcal I} r_\omega\bigr).
\end{equation*}
These identities show why the factor $\sigma^{-2}$ in $\partial_\varrho$ does not destroy the uniformity of the bounds. Indeed, the rescaled change of coordinates~$(x_1, x_2) \mapsto y/\sigma$, together with its inverse, is uniformly bounded at each fixed smooth order for positive~$\sigma$ close to zero. We can therefore express all these estimates in any fixed rescaled spatial norm. More precisely, we shall write $C^s_\sigma$ for the norm of the physical tensor components considered as functions of $(y/\sigma, Z)$. More precisely, for a scalar function $f$ we set
\[
\|f\|_{C^s_\sigma}
:= \|f(\sigma\,\cdot,Z)\|_{C^s},
\]
where the norm on the right is taken in the rescaled
coordinates $(y/\sigma,Z)$; for tensor fields, the
norm is defined componentwise in the physical coordinate
frame $(dy_1,dy_2,dZ)$ and its dual. This choice of norm is essential: the estimates do not give uniform bounds for arbitrarily many unscaled derivatives with respect to $\varrho$ as $\sigma \to 0$.

There is one feature of these coordinates that we must retain in the linear analysis: the leading metric differs from the product metric. Since $u = 2\sigma x_2 + O(\sigma^2)$, the coordinate differentials satisfy $dR = dy_1 + O(\sigma) \abs{dy}$, $dz = dy_2 + O(\sigma) \abs{dy}$, and $d\varphi = -dZ - 2\kappa_*\,dy_2 + O(\sigma) \abs{dy}$. It follows that, in any fixed~$C^s_\sigma$ norm,
\begin{equation}
G_0 = dy_1^2 + dy_2^2 + (dZ + c_* dy_2)^2 + O(\sigma), \qquad c_* = 2\kappa_*.
\label{kg24}
\end{equation}
As the shear $c_*$ is of order one, it cannot be discarded in this comparison. We shall keep it in the comparison metric for the coercivity argument of Section~\ref{kg41}.

We conclude this section by extending the coefficients beyond the physical domain. The reason for doing so is that the linear argument will require an inverse of curl without boundary conditions on the displacement. We shall obtain it on a closed manifold by identifying the transverse disk with a cap in a sphere. Specifically, choose $L > 1$, let $R_\Sigma = L\sigma$, and set
\begin{equation*}
M = \Sigma \times \T_N, \qquad \Sigma = S^2_{R_\Sigma}, \qquad \bar g = g_\Sigma + dZ^2.
\end{equation*}
The parameter $L$ compares the radius of the auxiliary sphere
with the transverse size of the reference torus. The latter
is of order $\sigma$, whereas the major radius of the reference
torus is of order one. Taking $L$ large makes the transverse
disk correspond to a small, nearly flat cap of $\Sigma$. We first fix such an $L$ and then take $\sigma$ sufficiently small that
$R_\Sigma=L\sigma$ is also small. Let $\vartheta$ and $\Theta$ be the polar and azimuthal angles on $\Sigma$. The action coordinate~$\varrho = R_\Sigma^2(1 - \cos\vartheta)$ puts the volume form in the form~$\mu = d\varrho \wedge d\Theta \wedge dZ$, so we can identify the physical source with the cap~$\{\varrho \leq \varrho_b\} \subset M$ through the coordinates~$(\varrho, \Theta, Z)$. Notice that rotation in $\Theta$ and translation in $Z$ leave this cap invariant. Consequently, restriction to the physical source commutes with averaging under the two actions. The sphere metric, expressed in these coordinates, is
\begin{equation*}
g_\Sigma = \frac{d\varrho^2}{2\varrho - \varrho^2/R_\Sigma^2} + \left(2\varrho - \frac{\varrho^2}{R_\Sigma^2}\right) d\Theta^2.
\end{equation*}
Denote by $h_\Sigma$ the second Cartesian coordinate of the embedded sphere. Its expression on the north cap is $h_\Sigma = y_2 \sqrt{1 - \abs{y}^2/(4R_\Sigma^2)}$. We can then introduce the comparison metric
\begin{equation}
g_* = g_\Sigma + (dZ + c_* dh_\Sigma)^2, \qquad c_* = 2\kappa_*,
\label{kg25}
\end{equation}
whose volume form is $\mu$ and which is uniformly equivalent to $\bar g$ for fixed~$\kappa_*$. Moreover, $\abs{\nabla_\Sigma h_\Sigma} \leq 1$, and on the cap the relative difference between $g_*$ and the principal metric in~\eqref{kg24} is $O(L^{-2})$. This choice thus keeps the shear of the physical metric and, at the same time, permits the use of the product Hodge operators. In what follows, the extension bounds will be expressed in $C^s_{R_\Sigma}$, namely the norm obtained after dilating the sphere to unit radius, using physical orthonormal tensor components and keeping $Z$ unscaled.

\begin{lemma}\label{kg26}
For fixed sufficiently large~$L$, the reference profiles and metric have smooth $U$-invariant extensions to $M$ satisfying
\begin{gather}
\omega = 1 + O(R_\Sigma^2 + \sigma^2), \qquad \omega' = 6 + O(\sigma), \qquad \beta_0 = 2 + O(R_\Sigma^2 + \sigma^2), \qquad \abs{\beta_0'} \leq C,
\label{kg27}
\\
(1 - \epsilon) g_* \leq G_0^{\rm ext} \leq (1 + \epsilon) g_*, \qquad \epsilon \leq C_{\kappa_*}(L^{-2} + \sigma).
\label{kg28}
\end{gather}
Their norms in rescaled coordinates are bounded uniformly in $\sigma$ at each fixed order. At fixed~$\sigma$, perturbations admit fixed smooth tame extension operators, equivariant under the two circle actions.
\end{lemma}

\begin{proof}
Let us apply Seeley's extension operator~\cite{Seeley} in a collar whose width is fixed in the rescaled cap. For the first profile, we extend $\omega' - 6$, multiply the extension by a cutoff, and integrate using its initial value on the physical domain. Since the interval is $0 \leq \varrho \leq 2R_\Sigma^2$, this procedure gives the first two estimates in the statement. Extending $\beta_0'$ and integrating in the same manner gives the remaining profile estimates. For the metric we instead extend the difference $G_0 - g_*$ and cut it off. Its relative $C^0$ bound is supplied by~\eqref{kg24}, so positivity is preserved with the choices specified in the lemma. The extension operator is bounded in every finite $C^s$ and Sobolev norm. To see that the rescaled estimates are uniform, notice that after dilation to the unit sphere the collars and cutoffs have sizes depending only on $L$. Finally, we average the extension operators over the two circle actions. This makes the operators equivariant while leaving their values on the cap unchanged, and completes the proof.
\end{proof}

\section{The nonlinear problem}
\label{kg29}

In this section we formulate the deformation problem on the physical source $\Omega_{\sigma, N}$. The magnetic field~$b$ and the volume form~$\mu$ will remain fixed; the unknown is a deformation of the covering immersion~$\Phi_0$. Following the discussion in Section~\ref{kg7}, we impose invariance under the two circle actions on the current of the transported field. We first record what this invariance implies for a divergence-free vector field on the source.

Given a tensor~$T$ on $\overline\Omega_{\sigma, N}$, we denote by $\av T$ the average obtained from the commuting flows of $X = \partial_\Theta$ and $U = \partial_Z$. Explicitly,
\begin{equation*}
\av T = \frac{1}{4\pi^2 N} \int_0^{2\pi} \int_0^{2\pi N}(\mathrm{Fl}_X^s \circ \mathrm{Fl}_U^t)^* T\,dt\,ds.
\end{equation*}
The range of this projection consists precisely of the tensors invariant under both flows. We shall also use the fact that averaging preserves divergence-free vector fields, which follows because both flows preserve $\mu$.

\begin{lemma}
\label{lem:profiles}
Every smooth divergence-free field that is invariant under $X$ and $U$ has the form
\begin{equation*}
j = \alpha(\varrho) b + \beta(\varrho) U,
\end{equation*}
with functions $\alpha$ and $\beta$ smooth up to the magnetic axis. Furthermore, for every integer $s\geq 0$,
\begin{align*}
\|\alpha\|_{H^s}
+\|\beta\|_{H^s}
&\leq C_s\|j\|_{H^{s+1}}\,,\\
\|\alpha'\|_{H^s}
+\|\beta'\|_{H^s}
&\leq C_s\|j\|_{H^{s+3}}\,.
\end{align*}
Here primes denote derivatives with respect to $\varrho$,
and all norms are taken on $\Omega_{\sigma,N}$, with the
functions evaluated at $\varrho=|y|^2/2$.
\end{lemma}

\begin{proof}
By invariance under $X$ and $U$, the horizontal component
of $j$ can be written, away from the axis, as
\[
j^y = a(\varrho)y + g(\varrho)Jy,
\]
while $j^Z$ depends only on $\varrho$. Since
$\mu=dy_1\wedge dy_2\wedge dZ$, the divergence equation
implies $(j^\varrho)'=0$. Moreover, $j^\varrho$ vanishes
at the axis by smoothness, so $j^\varrho\equiv0$.
As $j^\varrho=2\varrho a(\varrho)$, we conclude that
$a=0$. Recalling that
$b=\omega(\varrho)(X-\kappa_*U)$ and $X=Jy$, we obtain
\[
j=\alpha(\varrho)b+\beta(\varrho)U,
\qquad
\alpha=\frac{g}{\omega},
\qquad
\beta=j^Z+\kappa_*g.
\]

It remains to justify smoothness at the axis and the
estimates for extracting these profiles and their radial
derivatives. Consider a horizontal field
$v(y)=g(\varrho)Jy$ and a radial function
$H(y)=H_0(\varrho)$, where $\varrho=|y|^2/2$.
A direct computation gives
\[
\curl_2 v=2(g+\varrho g'),
\qquad
\Delta_y H=2(H_0'+\varrho H_0''),
\]
where $\curl_2 v=\partial_{y_1}v_2-\partial_{y_2}v_1$
and primes denote derivatives with respect to $\varrho$.
Consequently,
\[
\frac12(\curl_2 v)(\sqrt{t}\,y)
 =\frac{d}{dt}\bigl(tg(t\varrho)\bigr),
\qquad
\frac12(\Delta_y H)(\sqrt{t}\,y)
 =\frac{d}{dt}\bigl(tH_0'(t\varrho)\bigr).
\]
Integrating from $0$ to $1$, we obtain
\begin{equation}
g(\varrho(y))
 =\frac12\int_0^1(\curl_2 v)(\sqrt{t}\,y)\,dt,
\qquad
H_0'(\varrho(y))
 =\frac12\int_0^1(\Delta_y H)(\sqrt{t}\,y)\,dt.
\label{kg31}
\end{equation}

These formulas extend smoothly across $y=0$, since
their integrands are smooth functions of $y$.
The resulting functions are rotationally invariant,
and hence smooth functions of $\varrho$ up to the axis.

Finally, dilation by $\sqrt{t}$ costs at most
$t^{-1/2}$ in $L^2$ on the disk, while each spatial
derivative contributes an additional factor $\sqrt{t}$.
Since $t^{-1/2}$ is integrable on $(0,1)$, the formulas
above give Sobolev estimates with a loss of one
derivative for $g$ and two for $H_0'$.
Applying the second identity with $H=g$ yields a loss
of three derivatives for $g'$. Taking $H=j^Z$ and
recalling the expressions for $\alpha$ and $\beta$
completes the proof.
\end{proof}

Consider now an orientation-preserving immersion~$\Phi: \overline\Omega_{\sigma, N} \to \R^3$ sufficiently close to $\Phi_0$. We use the pullback metric~$G_\Phi = \Phi^* g_{\rm Eucl}$ and define the volume density by $\Phi^* dV = \rho_\Phi \mu$. For a one-form~$\eta$, the notation $\curl_\mu \eta$ will stand for the field determined by $\iota_{\curl_\mu \eta} \mu = d\eta$. Thus we can associate with the immersion the field
\begin{equation*}
c_\Phi = \curl_\mu\bigl(G_\Phi(b, \cdot)\bigr)
\end{equation*}
which is divergence-free with respect to $\mu$ for any immersion~$\Phi$, as follows from $d^2 = 0$. When $\rho_\Phi = 1$, this field is precisely the Euclidean current pulled back on an open set where $\Phi$ is invertible. We are therefore led to impose the equation
\begin{equation}
\mathcal F(\Phi) := \bigl(\rho_\Phi - 1, (\Id - \av) c_\Phi\bigr) = 0.
\label{eq:nonlinear-map}
\end{equation}
By~\eqref{kg17} and~\eqref{kg22}, the reference immersion is a solution of this equation. The usefulness of the invariance condition is that Lemma~\ref{lem:profiles} fixes the direction of the force, so that the pressure can be taken to depend only on $\varrho$. We make this implication precise in the following proposition.

\begin{proposition}
\label{kg33}
Let $\Phi$ be a smooth orientation-preserving immersion of $\overline\Omega_{\sigma, N}$ satisfying~\eqref{eq:nonlinear-map}. Then $c_\Phi = \alpha(\varrho) b + \beta(\varrho) U$ for smooth functions $\alpha, \beta$. For $0 < a \leq \varrho_b$, set
\begin{equation}
p(\varrho) = \int_\varrho^a \beta(s) \omega(s)\,ds.
\label{kg34}
\end{equation}
On each local inverse branch, the field~$B \circ \Phi = D\Phi\,b$ and pressure $P \circ \Phi = p(\varrho)$ satisfy $\Div B = 0$ and $(\curl B) \times B = \nabla P$. If $\Phi$ is an embedding on $\{\varrho \leq a\}$, these formulas define a smooth equilibrium on $\Omega = \Phi(\{\varrho < a\})$ satisfying~\eqref{kg2}, with $B$ nowhere zero and $P = 0$ on $\partial\Omega$. For $\Phi$ sufficiently close to $\Phi_0$, one has $p' < 0$. On an embedded sublevel the regular pressure levels are then nested tori surrounding the magnetic axis.
\end{proposition}
\begin{proof}
Since $c_\Phi$ is divergence-free by construction and satisfies $c_\Phi=\av c_\Phi$ by \eqref{eq:nonlinear-map}, the required expression for $c_\Phi$ follows from Lemma~\ref{lem:profiles}. Let us first verify the equilibrium equations on a local inverse branch of $\Phi$. As $\Phi^* dV = \mu$ and $\Div_\mu b = 0$, the transported field satisfies $\Div B = 0$. To identify its current, observe that $\Phi^*(B^\flat) = G_\Phi(b, \cdot)$. Exterior differentiation then yields $\Phi^*(\iota_{\curl B} dV) = \iota_{c_\Phi} \mu$, and hence $D\Phi\,c_\Phi = (\curl B) \circ \Phi$. It follows that the force satisfies
\begin{equation}
\Phi^*\bigl(((\curl B) \times B)^\flat\bigr) = \iota_b \iota_{c_\Phi} \mu = \beta\,\iota_b \iota_U \mu = -\beta\omega\,d\varrho = dp.
\label{kg35}
\end{equation}
To obtain the displayed identity we have used~\eqref{kg20} away from the axis. Both sides extend smoothly across the axis, so the force equation holds on the entire branch. If $\Phi$ is an embedding on the closed sublevel, these formulas define one pair~$(B, P)$ throughout its image. Since $b$ does not vanish and $d\varrho(b) = 0$, the transported field~$B$ is also nonvanishing and tangent to the boundary. The additive constant chosen for the pressure ensures that $P = 0$ there.

It remains to describe the pressure levels. At the reference immersion, $\beta = \beta_0 > 0$ by~\eqref{kg23}. Applying the profile extraction of Lemma~\ref{lem:profiles}, we see that $\beta$ remains positive for a solution sufficiently close to $\Phi_0$. Consequently $p' = -\beta\omega < 0$, and the pressure $P$ is positive in $\Omega$. The expression $\varrho = \abs{y}^2/2$ shows that its critical set consists only of the image of the axis and that its Hessian is negative definite in the normal directions there. All the other nonempty levels are the images of the tori $\varrho = \mathrm{constant}$, as claimed.
\end{proof}

As we know, the reference immersion satisfies $\mathcal{F}(\Phi_0)=0$ but it does not produce the desired knot. Our strategy will be to find a deformation direction that splits the magnetic axis into the desired knot and then use the Nash-Moser implicit function theorem to produce a family of immersions that both achieve the knot and satisfy $\mathcal{F}(\Phi)=0$. The first step is to study the linearization of the equation. Since $D\Phi$ is invertible, any variation of the immersion can be written as $\dot\Phi=D\Phi\,\xi$ for some vector field $\xi$ on the source torus. We impose no tangency condition on $\xi$ at $\partial \Omega_{\sigma, N}$ because the boundary of the image domain should be free to move.
\begin{proposition}
\label{kg36}
Let $\Phi$ be a smooth orientation-preserving immersion and $\dot\Phi = D\Phi\,\xi$. Then
\begin{equation*}
D\mathcal F(\Phi)[\dot\Phi] = \bigl(\Div_\mu(\rho_\Phi \xi), (\Id - \av) \dot c_\Phi\bigr),
\end{equation*}
where
\begin{equation*}
\dot c_\Phi = \curl_\mu\bigl(G_\Phi([b, \xi], \cdot)\bigr) - [c_\Phi, \xi] + (\Div_\mu \xi) c_\Phi.
\end{equation*}
\end{proposition}

\begin{proof}
In the interior we may realize the variation by composing $\Phi$ with the local flow of $\xi$; the derivative of this composition is $D\Phi\,\xi$. Differentiating the pullbacks of the Euclidean metric and of its volume form, we obtain $\dot G_\Phi = \mathcal L_\xi G_\Phi$ and $\dot\rho_\Phi = \Div_\mu(\rho_\Phi \xi)$. Smoothness extends these identities to the boundary. For the current, we must also use the fact that the field~$b$ is fixed, which gives
\begin{equation*}
\dot G_\Phi(b, \cdot) = \mathcal L_\xi\bigl(G_\Phi(b, \cdot)\bigr) + G_\Phi([b, \xi], \cdot).
\end{equation*}
As exterior differentiation commutes with $\mathcal L_\xi$, the exterior derivative of the first term on the right is
\begin{equation*}
\mathcal L_\xi(\iota_{c_\Phi} \mu) = \iota_{[\xi, c_\Phi]} \mu + (\Div_\mu \xi) \iota_{c_\Phi} \mu.
\end{equation*}
The second term contributes the curl appearing in the statement. Finally, the averaging operator~$\av$ is fixed on the source, so applying $\Id - \av$ gives the second component of the derivative and completes the proof.
\end{proof}

When $\mathcal F$ vanishes, we have $\rho_\Phi = 1$ and $c_\Phi = \bar j_\Phi$, with $\bar j_\Phi := \av c_\Phi$. We shall make these same two substitutions in the derivative at any nearby immersion, and write $M_\Phi$ for the operator so obtained. The advantage is that its current coefficient has the radial form of Lemma~\ref{lem:profiles}, whether or not the immersion solves the nonlinear equation. To keep track of the terms discarded by this substitution, Proposition~\ref{kg36} gives the exact formula
\begin{equation}
D\mathcal F(\Phi)[D\Phi\,\xi] = M_\Phi \xi + \err(\mathcal F(\Phi), \xi),
\label{kg37}
\end{equation}
where the terms discarded are collected in an error that is bilinear in the residual and the displacement:
\begin{equation}
\err((v, e), \xi) = \bigl(\Div_\mu(v\xi), (\Id - \av)(-[e, \xi] + (\Div_\mu \xi) e)\bigr).
\label{kg38}
\end{equation}
Indeed, the coefficient differences produced by the two substitutions are $v = \rho_\Phi - 1$ and $e = c_\Phi - \bar j_\Phi$. We shall use that, if $e$ is solenoidal, the second component of the error is solenoidal as well. This follows from the identity~$\Div_\mu[e, \xi] = e(\Div_\mu \xi)$.

We are now in a position to state the linear problem. Given a scalar function~$h$ and a vector field~$g$ on the source,  in Section~\ref{kg41} we shall solve for the displacement $\xi$ in the equation
\begin{equation}
M_\Phi \xi = \left(\Div_\mu \xi, (\Id - \av) \left\{\curl_\mu\bigl(G_\Phi([b, \xi], \cdot)\bigr) - [\bar j_\Phi, \xi] + (\Div_\mu \xi) \bar j_\Phi\right\}\right) = (h, g).
\label{kg39}
\end{equation}
The preceding divergence identity applies also to $\bar j_\Phi$. It shows that the second component of $M_\Phi \xi$ is solenoidal; by its definition this component has zero average as well. The appropriate fixed space for the data is therefore
\begin{equation}
\mathcal Y^\infty = C^\infty(\overline\Omega_{\sigma, N}) \times \{g \in C^\infty(T\overline\Omega_{\sigma, N}): \Div_\mu g = 0, \ \av g = 0\}.
\label{kg40}
\end{equation}
Our task is to construct a right inverse on the whole of this space without imposing a boundary condition on $\xi$. Once this has been done, the error identity~\eqref{kg37} will let us use that inverse to solve the nonlinear problem.

\section{Right inverse of the linear operator}
\label{kg41}
The goal of this section is to prove that the operator $M_\Phi$ defined in \eqref{kg39} has a right inverse that, when applied to smooth data, gives smooth fields with suitable dependence on the data and on $\Phi$.
\begin{proposition}
\label{prop right inverse of the linear operator}
Fix $\kappa_* > 1$ and $N \geq 2$. For every sufficiently small fixed band width $\delta > 0$ in~\eqref{kg46}, and all sufficiently small $\sigma>0$ depending on this choice, the operator~$M_\Phi$ has a smooth tame family of right inverses on a neighborhood of $\Phi_0$ specified by finitely many Sobolev norms,
\begin{equation*}
R_\Phi: \mathcal Y^\infty \longrightarrow C^\infty(T\overline\Omega_{\sigma, N}), \qquad M_\Phi R_\Phi = \Id.
\end{equation*}
More precisely, for some fixed integers $s_0, d$, the estimates have the form
\begin{equation}
\norm{R_\Phi h}_{H^s} \leq C_s\left(\norm{h}_{H^{s+d}} + (1 + \norm{\Phi - \Phi_0}_{H^{s+d}}) \norm{h}_{H^{s_0+d}}\right), \qquad s \geq s_0.
\label{kg43}
\end{equation}
The thickness $\sigma$ and the neighborhood of $\Phi_0$ are independent of $s$.
\end{proposition}

To prove the proposition, we first solve the divergence equation using a radial homotopy. Once we subtract a particular solution and absorb its contribution into the second datum, the equation we need to solve is
\begin{equation}
(\Id - \av) \{\curl_\mu(G([b, \xi], \cdot)) - [j, \xi]\} = g, \qquad \Div_\mu \xi = 0,
\label{reduced linear eq}
\end{equation}
where $G$ is close to the reference metric and $j = \alpha(\varrho) b + \beta(\varrho) U$ is close to $j_0$. The difficulty with a direct inversion of the magnetic transport $[b, \cdot]$ is that the divisors $n - \kappa_* m$ may vanish. To deal with these modes, we keep the radial displacement as an unknown on a fixed band containing the zeros of these divisors, and use transport along $j$ to recover $\xi$. This leads to a coercive variational problem whose energy norm also controls the commutators that we shall use to prove smooth regularity.

The slope~$\kappa_*$ and the covering degree~$N$ will remain fixed throughout this section. In particular, our constants are allowed to depend on $N$, although we will not always display this dependence.

The proof is organized as follows. We first derive the variational formulation and establish coercivity. Having obtained a weak solution, we prove the regularity estimates and reconstruct the physical displacement, which completes the proof of the proposition. We shall then extract a sharper estimate on the axis for a finite collection of longitudinal sectors; this is the estimate needed for the deformation of the axis in Section~\ref{kg90}.

\subsection{The variational problem}
To avoid imposing boundary conditions on $\xi$, we extend the circular cross-section of $\Omega_{\sigma, N}$ to a sphere $\Sigma$ and work on the closed manifold $M = \Sigma \times \T_N$, as discussed in Section~\ref{sec reference}. We use the extensions provided by Lemma~\ref{kg26} and extend the magnetic field as $b = \omega V$, where $V = X - \kappa_* U$. These profiles and the magnetic field are fixed throughout the construction. Meanwhile, we take the variable coefficients to be
\[p = (G, \alpha, \beta), \qquad j = \alpha(\varrho) b + \beta(\varrho) U,\]
in a sufficiently small finite-order neighborhood of $(G_0^{\rm ext}, 0, \beta_0)$, so in particular $\beta > 0$. The reader should note that, outside the physical cap, we only require $j$ to preserve volume and commute with $b$. There is no requirement that it be the current of an extended immersion.

Unless a different metric is explicitly indicated, the inner products and Sobolev norms below are defined using $\bar g$. We take inner products to be linear in their first argument, and use $\bar g$ to identify vector fields with one-forms via the musical isomorphisms. With these conventions, let us set
\[\HH = \{w \in L^2(TM \otimes \C): \Div_\mu w = 0, \ \av w = 0\}.\]
We denote by $\Pi$ the orthogonal projection onto $\HH$ and write $\CC = \curl_{\bar g}$. Since every harmonic one-form on $M$ is a multiple of $dZ$, vanishing-average condition removes this harmonic part. We can therefore apply the Hodge decomposition to obtain
\begin{equation}
\CC^{-1} = \CC\Delta_1^{-1} \quad \hbox{on }\HH, \qquad \CC^{-1}(\Id - \av) \curl_\mu \eta = \Pi\eta^\sharp.
\label{kg45}
\end{equation}
Here we define the inverse Laplacian to vanish on harmonic forms. In the second identity, we have used that $\curl_\mu$ and $\CC$ are defined with the same volume form. On general one-forms we extend the inverse by $\CC^{-1} \Pi = \CC\Delta_1^{-1}(\Id - \av)$, and keep the notation $\CC^{-1}$ for this extension. 

With these definitions, taking $\CC^{-1}$ in \eqref{reduced linear eq} and introducing a new variable $w\coloneqq [b,\xi]$, we obtain
\[\Pi\, G^\sharp w-\CC^{-1}\mathcal{L}_j\xi=\tilde{g},\]
where we have also modified the datum on the right-hand side. Note that the first operator on the left-hand side is positive. If we could obtain $\xi$ from $w$, we would then have an equation only in terms of $w$, and we could hope to get coercivity from the positive energy given by the first operator. Unfortunately, at exact resonances we cannot invert $\mathcal{L}_b$ to recover $\xi$. Furthermore, close to resonances the second operator can be too large to be controlled by the first. 

The key will be to keep the radial component of $\xi$ as an independent unknown, instead of trying to obtain it from $w$. Regarding the second operator, we will write it as the sum of a positive operator (which may be large) and a negative operator that can be controlled by making use of the divergence-free condition. However, we first consider the simpler case and work far from resonances.

For a tensor $T$ we define the angular projections
\[\mathsf{P}_{nm}T\coloneqq \frac{1}{4\pi^2N}\int_0^{2\pi}\int_0^{2\pi N}e^{-\ii (ns+mt)}\left(\rm{Fl}_X^s\circ \rm{Fl}_U^t\right)^\ast\hspace{1pt}T\,dt\,ds.\]
In particular, we have $\mathcal{A}=\mathsf{P}_{00}$. Since the flows of $X$ and $U$ commute, we have:
\[\mathcal{L}_X(\mathsf{P}_{nm}T)=\ii n\mathsf{P}_{nm}T, \qquad \quad \mathcal{L}_U(\mathsf{P}_{nm}T)=\ii m\mathsf{P}_{nm}T \qquad \forall n\in\Z, m\in N^{-1}\Z.\]
Hence, on a mode $(n,m)$, transport by $V=X-\kappa_* U$ reduces to:
\begin{equation}
\label{transport by V}
\mathcal{L}_V(\mathsf{P}_{nm}T)=\ii (n-\kappa_* m)\mathsf{P}_{nm}T\equiv \ii k \mathsf{P}_{nm}T.
\end{equation}
Thus, inverting $\mathcal{L}_V$ will be problematic when $k\coloneqq n-\kappa_* m$ is too small. To remove these modes, we fix a positive number $2\delta < \min(1, \kappa_*/N)$ that will be chosen later and we define
\begin{equation}
\mathsf{P}\coloneqq\sum_{(n,m)\neq (0,0),\ \abs{k}<\delta} \mathsf{P}_{nm}, \hspace{70pt} \mathsf E \coloneqq \Id - \av - \mathsf P.
\label{kg46}
\end{equation}
\begin{remark}
\label{rem width band}
The width $\delta>0$ has been chosen small enough so that for modes in the band $\mathsf P$ both $n$ and $m$ are nonzero, there is at most one $n$ for each $m$, and $\abs{m/n} \leq C_{\kappa_*}$ for some $C_{\kappa_*}>0$.  
\end{remark}
We choose $\delta$ before fixing $\sigma$. More precisely, the band is taken narrow enough that it contains only the exact resonances in the finitely many longitudinal sectors that will be used in Section~\ref{kg90}. We do not impose any lower bound on the nonzero frequencies that remain inside the band.

It follows from \eqref{transport by V} that, away from the band, transport by $V$ has the bounded inverse
\begin{equation*}
Q_E \coloneqq \sum_{\abs{k}\geq\delta} \frac{\mathsf P_{nm}}{\ii k}.
\end{equation*}
Since this operator commutes with the product Hodge Laplacians, the same bound $\delta^{-1}$ holds for its norm at every Sobolev order. To invert transport by $b=\omega V$, we also have to take into account the variation of $\omega$:
\begin{lemma}\label{kg47}
For $e \in \mathsf E\HH$, let
\begin{equation}
T_b e \coloneqq \omega^{-1} Q_E e + \frac{\omega'}{\omega^2} Q_E^2(e^\varrho) V.
\label{kg48}
\end{equation}
Then $T_b e$ is divergence-free and $[b, T_b e] = e$. The operator~$A_E = \CC^{-1} \mathcal{L}_j T_b$ on $\mathsf E\HH$ is bounded and self-adjoint, with
\begin{equation}
\norm{A_E}_{L^2\to L^2} \leq C_{\kappa_*} \{R_\Sigma(1 + \delta^{-1}) + R_\Sigma^2(\delta^{-1} + \delta^{-2})\} \leq C_{\kappa_*, \delta} R_\Sigma
\label{kg49}
\end{equation}
for $R_\Sigma \leq 1$ and coefficients in fixed low bounds.
\end{lemma}
\begin{proof}
The first claim follows by using the identity~$[b, \xi] = \omega[V, \xi] - \omega'\xi^\varrho V$. Indeed, the two terms of~\eqref{kg48} have divergences $-\omega'\omega^{-2} Q_E(e^\varrho)$ and its opposite, so their sum is divergence-free. To prove symmetry, it is convenient to work with the helicity pairing $\mathsf{h}(u,v) = \langle\CC^{-1} u, v\rangle$ on exact divergence-free fields. Cartan's magic formula, followed by integration by parts, leads to
\[\mathsf{h}([Y,u], v) + \mathsf{h}(u,[Y,v]) =0\]
for any real solenoidal field~$Y$. Since $[b, j] = 0$, we see that $T_b$ commutes with $\mathcal{L}_j$. Both operators are skew for the helicity pairing, and their product is therefore symmetric, as we wanted to prove.

Next, we estimate the inverse curl on a nonzero weight. For a divergence-free field, the relevant sphere eigenvalues are bounded below by $c(1 + n^2)/R_\Sigma^2$. When $n = 0$ and $m \neq 0$, one must also rule out a sphere-constant vertical component. This is the only possible exception, and it is excluded by integrating the divergence equation over the sphere. We thus obtain
\begin{equation}
\bigl\|\CC^{-1} \mathsf P_{nm}\bigr\|_{L^2\to L^2} \leq C\bigl((1 + n^2)/R_\Sigma^2 + m^2\bigr)^{-1/2}.
\label{kg50}
\end{equation}
A direct expansion gives
\begin{equation}
\mathcal{L}_j T_b = \alpha \mathsf E + \frac\beta\omega UQ_E + \frac{\beta\omega'}{\omega^2}(V \otimes d\varrho) UQ_E^2  -\frac{\alpha'} \omega(b \otimes d\varrho) Q_E - \frac{\beta'} \omega(U \otimes d\varrho) Q_E,
\label{kg51}
\end{equation}
where $U$, when applied to a field, denotes Lie differentiation. Using $\kappa_* m = n - k$ in~\eqref{kg50}, we infer that $\abs{m} \bigl\|\CC^{-1} \mathsf P_{nm}\bigr\|_{L^2\to L^2} \leq C_{\kappa_*} R_\Sigma(1 + \abs{k})$. As $\abs{d\varrho} \leq R_\Sigma$ and $V$ is bounded, this estimate applied to the terms of the expansion proves~\eqref{kg49} from~\eqref{kg51}. Here the coefficients are radial and hence preserve the weights, so the modal estimates can be summed orthogonally. The boundedness just established also allows us to extend symmetry from smooth fields to $\mathsf E\HH$.
\end{proof}

We now turn to the modes selected by the band, where it is preferable to keep the radial displacement as an unknown. Given $w = [b, \xi]$, let us introduce $f = \omega\xi^\varrho$ and $\tau = w^Z$. The solenoidal condition then determines the components of the field as
\begin{equation}
w^\varrho = \ii kf, \qquad w^\Theta = -\frac{k}{n} f' - \frac mn\tau, \qquad w^Z = \tau.
\label{kg52}
\end{equation}
Although written in angular coordinates, these formulas also make sense globally at the poles. To express them in an intrinsic form valid on all $M$, 
let $\mu_\Sigma=d\varrho\wedge d\Theta$ denote the area form
of $(\Sigma,g_\Sigma)$, where $d_\Sigma$ is the differential on the sphere $\Sigma$. For a scalar function $f$ on
$\Sigma$, we denote by $X_f^\Sigma$ the associated
Hamiltonian vector field, defined by
$\iota_{X_f^\Sigma}\mu_\Sigma=d_\Sigma f$. For smooth $z=(e,f,\tau)$ we then define
\begin{equation}
J_1 f = \frac{k}{n} X_f^\Sigma, \qquad J_2 \tau = \tau U - \frac mn\tau X, \qquad  Sz = e + J_1 f + J_2 \tau.
\label{kg53}
\end{equation}
Both $J_1 f$ and $J_2 \tau$ are smooth divergence-free fields, so $Sz$ is smooth and divergence-free. Away from the poles, $J_1 f+J_2\tau$ gives exactly the components of $w$ in \eqref{kg52} for each mode $(n,m)$ in the band $\mathsf{P}$. To see this, note that in coordinates $(\varrho, \Theta)$ we have $X_f^\Sigma=(\partial_\Theta f)\partial_\varrho-(\partial_\varrho f)\partial_\Theta$. The purpose of the term $e$ is to supply the modes of $w$ outside the resonance band.

Notice that when $k = 0$, we have $J_1 = 0$ and $J_2 \tau = -(\tau/\kappa_*) V$. Keeping $f$ as an independent unknown allows us to recover the radial displacement even at a magnetic resonance. We now introduce the spaces where we will work:
\begin{lemma}
Let 
\begin{align}
\ZZ_0 &= \mathsf E\HH \oplus \mathsf PL^2(M) \oplus \mathsf PL^2(M), \\[2pt]
\ZZ &= \{z = (e, f, \tau) \in \ZZ_0: J_1 f \in L^2\}.
\label{kg54}
\end{align}
and equip $\ZZ$ with the norm 
$E_0(z)^2 = \norm{e}_0^2 + \norm{f}_0^2 + \norm{\tau}_0^2 + \norm{J_1 f}_0^2$. Then $\ZZ$ is a Hilbert space and the norm $E_0$ satisfies
\begin{equation}
E_0(z)^2 \sim \norm{Sz}_0^2 + \norm{f}_0^2.
\label{kg55}
\end{equation}
\end{lemma}
\begin{proof}
Since $J_1$ is closed on its maximal distributional domain and $\ZZ_0$ is a Hilbert space, so is $\ZZ$. To prove the equivalence of the norms, first observe that $e$ only contains modes outside the band $\mathsf{P}$ while $J_1f+J_2\tau$ only contains modes in the band, so:
\begin{equation}
\mathsf{E}Sz=e, \hspace{40pt} \mathsf{P}Sz=J_1f+J_2\tau.
\label{decomp Sz}
\end{equation}
Since these projections are orthogonal in $L^2$, we have
\begin{equation}
\|Sz\|_0^2=\|e\|_0^2+\|\mathsf{P}Sz\|_0^2
\label{orthogonality Sz}
\end{equation}
and, in particular, each term on the right is bounded by $\|Sz\|_0^2$. The next observation is that $J_1 f$ is tangent to the sphere, while $J_2\tau=\tau U-\frac{m}{n} X$ has $Z$-component $\tau$. Therefore, we can recover $\tau$ as:
\[\tau = dZ(\mathsf{P}Sz)\qquad \Rightarrow \qquad \|\tau\|_0\leq \|\mathsf{P}Sz\|_0\leq \|Sz\|_0.\]
Finally, since $\frac{m}{n}$ is uniformly bounded on the resonance band (see Remark~\ref{rem width band}), we see that $J_2$ is bounded on $L^2$. Hence,
\[\|J_1f\|_0=\|\mathsf{P}Sz-J_2\tau\|_0\leq \|\mathsf{P}Sz\|_0+\|J_2\tau\|_0\lesssim \|Sz\|_0.\]
Collecting all these bounds, we conclude $E_0(z)^2\lesssim \|Sz\|_0^2+\|f\|_0^2$. Note that the additional term $\|f\|_0^2$ is necessary because on an exact resonance $k=0$ we have $J_1\equiv 0$, so $Sz$ carries no information about $f$. For the reverse inequality, it follows from \eqref{orthogonality Sz} that it suffices to prove $\|\mathsf{P}Sz\|^2_0\lesssim\|\tau\|_0^2+\|J_1f\|_0^2$. This follows from \eqref{decomp Sz} and boundedness of $J_2$.
\end{proof}

To determine the coupling between $f$ and $\tau$, we use the current transport. Writing $j^\Theta = \alpha\omega$ and $j^Z = \beta - \kappa_* \alpha\omega$, we introduce the coefficients
\begin{equation}
\begin{gathered}
\lambda = nj^\Theta + mj^Z = \alpha\omega k + \beta m, \qquad \mathsf B = -\frac{\lambda}{n\omega}, \\
\mathsf C = \frac{k(j^Z)' + \kappa_* \lambda\omega'/\omega}{n\omega} = \frac{\beta\omega'}{\omega^2} + \frac kn\left[\left(\frac\beta\omega\right)' - \kappa_* \alpha'\right].
\end{gathered}
\label{kg56}
\end{equation}
These coefficients are bounded real radial modal multipliers. By decreasing the coefficient neighborhood if necessary, we can also ensure that
\begin{equation}
\abs{\lambda} \geq c \abs{m} \quad \hbox{on }\mathsf P.
\label{kg57}
\end{equation}
This follows because $\abs{m} \geq N^{-1}$, $\abs{k} < \delta$, and $\beta$ remains positive. Hence the current transport is invertible throughout the selected band. If we write $j_\varrho = (j^\Theta)'X + (j^Z)'U$, its inverse takes the form
\begin{equation}
T_j \mathfrak r \coloneqq \frac{\mathfrak r}{\ii\lambda} + \frac{\mathfrak r^\varrho}{(\ii\lambda)^2} j_\varrho.
\label{kg58}
\end{equation}
Indeed, the identity~$[j, \xi] = \ii\lambda\xi - \xi^\varrho j_\varrho$ immediately gives $[j, T_j \mathfrak r] = \mathfrak r$. When $\mathfrak r$ is solenoidal, we can take the divergence of this equation and use~\eqref{kg57} to conclude that $\Div T_j \mathfrak r = 0$ as well.

We now consider the Hermitian form
\begin{equation}
\cform_p(z, \widehat z) =\int_M G(Sz, \overline{S\widehat z})\,\mu - \langle A_E e, \widehat e\rangle + \langle\mathsf Cf, \widehat f\rangle +\langle\mathsf B\tau, \widehat f\rangle + \langle\mathsf Bf, \widehat\tau\rangle.
\label{kg59}
\end{equation}
This form is continuous on $\ZZ$. Assuming that there exists a smooth solution to the variational problem that it poses, we now show that it produces a solution to \eqref{reduced linear eq}, as desired. Existence of a solution in $\ZZ$ will be given by the Lax--Milgram theorem after we prove Proposition~\ref{kg69}. Smoothness will be proved in the next subsection.
\begin{lemma}\label{kg60}
Suppose $g \in \HH$ and a smooth $z = (e, f, \tau)$ satisfies
\begin{equation}
\cform_p(z, \widehat z) = \langle\CC^{-1} g, S\widehat z\rangle \qquad \quad\forall \widehat z \in \ZZ.
\label{kg61}
\end{equation}
Set
\begin{equation}
w = Sz, \qquad \mathfrak r = (\Id - \av) \curl_\mu(G(w, \cdot)) - g, \qquad \xi = T_b e + T_j \mathsf P\mathfrak r.
\label{kg62}
\end{equation}
Then $\xi$ is solenoidal, jointly mean-free, and satisfies
\begin{equation}
[b, \xi] = w, \qquad [j, \xi] = \mathfrak r.
\label{kg63}
\end{equation}
In particular it solves~\eqref{reduced linear eq}.
\end{lemma}

\begin{proof}
Taking the $\widehat e$ component as a test and using~\eqref{kg45}, we get $\mathsf E\CC^{-1} \mathfrak r = A_E e$, from which it follows that $\mathsf E\mathfrak r = [j, T_b e]$. It remains to understand the selected weights. On each such weight, the field~$J_1 \widehat f + J_2 \widehat\tau$ has the smooth potential
\begin{equation*}
\eta = \frac{\ii\widehat\tau}{n}\,d\varrho + \frac{k\widehat f}{n}\,dZ, \qquad d\eta = \iota_{J_1 \widehat f+J_2 \widehat\tau} \mu.
\end{equation*}
Using this potential in the pairing, we find that the tests in the two scalar components are exactly
\begin{equation}
\frac kn\mathfrak r^Z + \mathsf Cf + \mathsf B\tau = 0, \qquad - \frac\ii n\mathfrak r^\varrho + \mathsf Bf = 0.
\label{kg64}
\end{equation}
Let us set~$\xi_P = T_j \mathsf P\mathfrak r$. By the second equation, $\mathfrak r^\varrho = \ii\lambda f/\omega$, which gives $\xi_P^\varrho = f/\omega$ and hence $[b, \xi_P]^\varrho = \ii kf$. Once this radial component has been determined, the first equation yields
\begin{equation*}
[b, \xi_P]^Z = \frac{\omega k}{\lambda} \mathfrak r^Z + \frac{k(j^Z)'} \lambda f + \frac{\kappa_* \omega'} \omega f = \tau.
\end{equation*}
The preceding identities are valid also when $k = 0$. Thus $[b, \xi_P]$ and $\mathsf Pw$ can only differ in their poloidal component; since both fields are solenoidal and the mode satisfies $n \neq 0$, that component must vanish as well. Together with the off-band equation, this proves~\eqref{kg63}. Finally, the invariance of $j$ gives $\av[j, \xi] = [j, \av\xi] = 0$.
\end{proof}

Notice that at exact resonance the first equation in~\eqref{kg64} becomes $\mathsf Cf + \mathsf B\tau = 0$, so it determines $f = \omega\tau/(\kappa_* \omega')$ without requiring division by $k$. 

In the previous lemma we assume the existence of a smooth solution of \eqref{kg61}. We will now prove coercivity of $\mathfrak{b}_p$, which gives a solution in $\ZZ$. Smoothness will be proved in the next subsection. We begin with two technical lemmas that will allow us to make use of the divergence-free condition. For this, we denote by $\Pi_\Sigma$ the orthogonal projection onto divergence-free tangent fields on the sphere $\Sigma$.
\begin{lemma}\label{kg65}
If $\phi$ is a scalar function of poloidal weight~$n \neq 0$, its Hodge decomposition $\phi\nabla_\Sigma h_\Sigma = u_\Sigma + \nabla_\Sigma p_\Sigma$, with $u_\Sigma = \Pi_\Sigma(\phi\nabla_\Sigma h_\Sigma)$ and $\int_\Sigma p_\Sigma = 0$, satisfies
\begin{equation}
\norm{u_\Sigma}_0^2 \leq \tfrac12 \norm{\phi}_0^2, \qquad \norm{p_\Sigma}_0 \leq C\frac{R_\Sigma}{\abs{n}} \norm{\phi}_0.
\label{kg66}
\end{equation}
\end{lemma}

\begin{proof}
To prove the first estimate, let $\mathcal R_1, \mathcal R_2, \mathcal R_3$ denote the dimensionless rotation fields about the three Cartesian axes, chosen so that $\mathcal R_3 = X$. The scalar curl of $\phi\nabla_\Sigma h_\Sigma$ is, up to a sign, $R_\Sigma^{-1} \mathcal R_2 \phi$. Since there are no harmonic one-forms on the sphere, we can therefore write
\begin{equation*}
\norm{u_\Sigma}_0^2 = R_\Sigma^{-2} \langle\mathcal R_2 \phi, \Delta_\Sigma^{-1} \mathcal R_2 \phi\rangle.
\end{equation*}
The rotations commute with $\Delta_\Sigma$, and the identity $-\sum_{i=1}^3 \mathcal R_i^2 = R_\Sigma^2 \Delta_\Sigma$ relates their quadratic expressions. Now a rotation through $\pi/2$ about the third axis exchanges the first two terms while changing $\phi$ only by a constant of modulus one. These two terms must consequently be equal, and their sum is bounded by $R_\Sigma^2 \norm{\phi}_0^2$, which proves the first estimate. For the second one, observe that the potential $p_\Sigma$ only has the weights~$n - 1, n + 1$ and satisfies $\norm{\nabla_\Sigma p_\Sigma}_0 \leq \norm{\phi}_0$. When $\abs{n} \geq 2$, the nonzero sphere eigenvalues in these weights are at least $cn^2/R_\Sigma^2$. When $\abs{n} = 1$, the mean-zero normalization removes the constant eigenfunction and leaves the spectral gap $2/R_\Sigma^2$. The desired bound for the potential follows in both cases.
\end{proof}

\begin{lemma}\label{kg67}
For every solenoidal field~$w = v + w^Z U$ on $M$, with $v$ tangent to the sphere,
\begin{equation}
\norm{\mathsf Pw^Z}_0^2 \leq (1 + c_*^2/2 + C_{\kappa_*} R_\Sigma) \int_M g_*(w, \overline w)\,\mu.
\label{kg68}
\end{equation}
\end{lemma}

\begin{proof}
We prove the estimate in a fixed longitudinal sector $m$ for which the band selects a poloidal weight~$n$. With the notation $\phi = \mathsf P_{nm} w^Z$ and $\mathcal E_m = \norm{v}_0^2 + \norm{w^Z + c_* dh_\Sigma(v)}_0^2$, the decomposition in the preceding lemma gives
\begin{equation*}
\norm{\phi}_0^2 = \langle w^Z + c_* dh_\Sigma(v), \phi\rangle - c_* \langle v, u_\Sigma\rangle - c_* \langle v, \nabla_\Sigma p_\Sigma\rangle.
\end{equation*}
The sum of the first two terms is bounded by $\sqrt{1 + c_*^2/2}\,\mathcal E_m^{1/2} \norm{\phi}_0$. To estimate the last term, we integrate by parts and use $\Div_\Sigma v = -\ii mw^Z$, obtaining the bound $C \abs{c_*} R_\Sigma \abs{m/n} \norm{w^Z}_0 \norm{\phi}_0$. Since $\abs{m/n} \leq C_{\kappa_*}$ and $\norm{w^Z}_0 \leq \sqrt{1 + c_*^2}\,\mathcal E_m^{1/2}$, we can now square the resulting inequality and sum over $m$ to prove~\eqref{kg68}. In this last step we use orthogonality only in the longitudinal variable. All poloidal weights of $w$ are kept in the energy, as is necessary because of the shear in $g_*$.
\end{proof}

\begin{proposition}\label{kg69}
For each fixed~$\kappa_* > 1$ and $N \geq 2$, choose $L$ sufficiently large. For every sufficiently small fixed~$\delta > 0$ and then all sufficiently small~$\sigma > 0$, we have
\begin{equation}
\cform_p(z, z) \geq c_{\kappa_*} E_0(z)^2 \qquad (z \in \ZZ)
\label{kg70}
\end{equation}
throughout one sufficiently small finite-order coefficient neighborhood.
\end{proposition}

\begin{proof}
By~\eqref{kg27} and~\eqref{kg56}, on the band we have
\begin{equation*}
\mathsf C \geq c > 0, \qquad \abs{\mathsf B/\mathsf C} \leq C_{\kappa_*}, \qquad \frac{\mathsf B^2}{\mathsf C} \leq \frac1{3\kappa_*^2} + C_{\kappa_*}(\delta + \sigma + R_\Sigma^2 + \varepsilon_p),
\end{equation*}
where $\varepsilon_p$ measures the perturbation in a fixed low norm. The role of these inequalities becomes apparent upon completing the square in the scalar terms:
\begin{equation*}
\langle\mathsf Cf, f\rangle + 2\operatorname{Re}\langle\mathsf B\tau, f\rangle = \norm{\mathsf C^{1/2}(f + \mathsf C^{-1} \mathsf B\tau)}_0^2 - \langle\mathsf B^2 \mathsf C^{-1} \tau, \tau\rangle.
\end{equation*}
We have $\mathsf P(Sz)^Z = \tau$, so Lemma~\ref{kg67} bounds the negative term by the full metric energy. Combining this with metric comparison and~\eqref{kg49}, we find that the coefficient left in front of the metric energy is
\begin{equation}
1 - \frac{1 + 2\kappa_*^2}{3\kappa_*^2} = \frac{\kappa_*^2 - 1}{3\kappa_*^2} > 0,
\label{kg71}
\end{equation}
up to an error bounded by $C_{\kappa_*}(L^{-2} + \delta + \sigma + R_\Sigma^2 + \varepsilon_p) + C_{\kappa_*, \delta} R_\Sigma$. We choose the parameters in the order stated in the proposition and then shrink the coefficient neighborhood. With these choices, the positive terms that remain control
\begin{equation*}
\norm{Sz}_0^2 + \norm{f + \mathsf C^{-1} \mathsf B\tau}_0^2.
\end{equation*}
Since $\norm{\tau}_0 \leq \norm{Sz}_0$ and $\mathsf B/\mathsf C$ is bounded, this also controls $\norm{Sz}_0^2 + \norm{f}_0^2$. By~\eqref{kg55}, we have therefore obtained control of $E_0(z)^2$, which proves coercivity.
\end{proof}

We can now apply the Lax--Milgram theorem to obtain a unique weak solution for each continuous functional on $\ZZ$. Observe that the argument treats an exact resonance and a small nonzero magnetic frequency in the same way. Our next task is to show that this weak solution is smooth whenever the data are smooth, without changing the coefficient neighborhood.

\subsection{Regularity and proof of Proposition~\ref{prop right inverse of the linear operator}}

For the nonlinear argument, the weak solution just constructed must give a smooth tame inverse. We shall obtain the required estimates by commuting the variational equation with fixed Sobolev gradings. The point is that these commutators have lower order in the unknown. Their contribution can therefore be absorbed by interpolation, so the coefficient neighborhood chosen for coercivity will work at every higher order as well.

We use the scalar and one-form Hodge Laplacians $\Delta_0, \Delta_1$ of the product metric to define the Sobolev gradings
\begin{equation*}
\Lambda_0 = (1 + \Delta_0)^{1/2}, \qquad \Lambda_1 = (1 + \Delta_1)^{1/2}, \qquad \Lambda = \operatorname{diag}(\Lambda_1, \Lambda_0, \Lambda_0).
\end{equation*}
Here $\Lambda_1$ acts on the full bundle of one-forms. This is necessary because the fields to which it will be applied, such as $G^\sharp Sz$, need not be solenoidal. Up to the fixed musical identification, the sphere Hamiltonian map is $*_\Sigma d_\Sigma$, and hence it intertwines the gradings as follows:
\begin{equation}
\Lambda_1^s J_1 = J_1 \Lambda_0^s.
\label{kg72}
\end{equation}
The angular projections and modal multipliers commute with the corresponding gradings as well. For each $s \geq 0$, let us write
\begin{equation}
E_s(z)^2 = \norm{e}_s^2 + \norm{f}_s^2 + \norm{\tau}_s^2 + \norm{J_1 f}_s^2.
\label{kg73}
\end{equation}
With this notation, $E_s(z) = E_0(\Lambda^s z)$. Notice also that spectral truncations preserve all the constraints and are dense for each of these graph norms.

We will use the more general weak equation
\begin{equation}
\cform_p(z, \widehat z) = \langle H_1, S\widehat z\rangle + \langle H_2, \widehat z\rangle, \qquad \widehat z \in \ZZ,
\label{kg74}
\end{equation}
where $H_1$ is a vector field and $H_2 \in \ZZ_0$. We include the additional datum~$H_2$ because it will be needed when differentiating the equation with respect to the coefficients. For brevity, set~$D_s = \norm{H_1}_s + \norm{H_2}_s$. In the coefficient norm~$\norm{p}_s$ we include the lifted radial derivatives $\alpha', \beta'$ appearing in the form. By Lemma~\ref{lem:profiles}, this amounts to a fixed shift in the number of spatial derivatives.

\begin{proposition}\label{kg75}
There are fixed integers $s_0, d_0$ such that the solution of~\eqref{kg74} satisfies, for all $s \geq s_0$,
\begin{equation}
E_s(z) \leq C_s\bigl(D_s + (1 + \norm{p}_{s+d_0}) D_{s_0}\bigr).
\label{kg76}
\end{equation}
It is smooth for smooth data, and depends smoothly tamely on the coefficients and the data. All these assertions hold on one fixed finite-order coefficient neighborhood.
\end{proposition}

\begin{proof}
We first establish the energy estimate on finite-dimensional approximations. In this way the commutation of derivatives will be justified before any regularity of the weak solution has been proved. Let us write the form as
\begin{equation*}
\cform_p(z, \widehat z) = \langle G^\sharp Sz, S\widehat z\rangle + \langle\mathsf R_p z, \widehat z\rangle,
\end{equation*}
where $\mathsf R_p$ is formed by the bounded blocks $-A_E, \mathsf B, \mathsf C$ in~\eqref{kg59}. For ease of notation, in the formulas below we omit the musical identification on $G$.

We shall use the standard tame product and commutator
estimates in Sobolev spaces. Since $M$ is
three-dimensional, we fix an integer $s_0 \geq 2$,
so that $H^{s_0}$ is an algebra and
$H^{s_0+d_0}$ embeds into $W^{1,\infty}$ (we can take $d_0=1$ in what follows).
The estimates then read
\begin{align}
\norm{uv}_s &\leq C_s(\norm{u}_s \norm{v}_{s_0} + \norm{u}_{s_0} \norm{v}_s),
\label{kg77}
\\
\norm{[\Lambda^s, M_u] v}_0 &\leq C_s(\norm{u}_{s_0+d_0} \norm{v}_{s-1} + \norm{u}_{s+d_0} \norm{v}_{s_0}).
\label{kg78}
\end{align}
The scalar versions of these inequalities can be found in~\cite[(3.1.59) and Proposition~3.6.A]{Taylor}. To obtain the tensor versions, one works in local bundle frames. Since the Hodge gradings have scalar principal symbols, the additional matrix-valued terms have order at most $s - 1$ and are bounded by the product estimate. Sobolev embedding then gives the displayed inequalities.

By~\eqref{kg72}, the commutator $\mathcal K_s = \Lambda_1^s S - S\Lambda^s$ receives a contribution only from the $J_2$ block. Applying the commutator estimate to multiplication by the fixed smooth field~$X$, we therefore obtain
\begin{equation}
\norm{\mathcal K_s z}_0 \leq C_s' \norm{\tau}_{s-1}, \qquad \norm{\mathcal K_s^* v}_0 \leq C_s' \norm{v}_{s-1}.
\label{kg79}
\end{equation}
To estimate commutators with $\mathsf R_p$, we use formula~\eqref{kg51}. It gives the same bounds as a bounded operator obtained by composing smooth multiplications and fixed spectral multipliers. Indeed, the factors $U$ that appear to be first order can be moved past the $U$-invariant coefficients and placed next to $\CC^{-1}$. The resulting operator~$\CC^{-1} U$ is bounded on every Sobolev space. The other factors in~\eqref{kg51} and~\eqref{kg56} have order zero, and their fixed spectral parts commute with the gradings. It follows that
\begin{equation}
\norm{[\Lambda^s, \mathsf R_p] z}_0 \leq C_s\bigl(\norm{z}_{s-1} + (1 + \norm{p}_{s+d_0}) \norm{z}_{s_0}\bigr)
\label{kg80}
\end{equation}
on a fixed low-norm ball of coefficients. At the finitely many orders up to $s_0$, bounding finitely many derivatives of the coefficients instead gives the cruder estimate~$C \norm{z}_{s-1}$. The same argument applies after differentiation of the operators with respect to their coefficients.

We now justify the energy calculation by taking an increasing sequence of finite spectral subspaces $V_j \subset \ZZ$ and solving~\eqref{kg74} with $z_j, \widehat z \in V_j$. Coercivity gives $E_0(z_j) \leq CD_0$. Since $\Lambda^{2s} z_j$ also belongs to $V_j$, it is an admissible test in this finite-dimensional equation. The calculation is based on the following exact identity, valid for a smooth $z$:
\begin{align}
& \cform_p(z, \Lambda^{2s} z) - \cform_p(\Lambda^s z, \Lambda^s z) \notag \\
& \quad = \langle[\Lambda_1^s, G] Sz + G\mathcal K_s z, S\Lambda^s z\rangle - \langle\mathcal K_s^* GSz, \Lambda^s z\rangle + \langle[\Lambda^s, \mathsf R_p] z, \Lambda^s z\rangle.
\label{kg81}
\end{align}
We use this identity directly for the form, so we never have to differentiate the unbounded operator represented by it. In particular, the factor $J_1 f$ is estimated in its graph norm throughout. Applying~\eqref{kg77}--\eqref{kg80}, we bound the absolute value of the right-hand side by
\begin{equation*}
C_s\bigl(E_{s-1}(z) + (1 + \norm{p}_{s+d_0}) E_{s_0}(z)\bigr) E_s(z).
\end{equation*}
For the right-hand side of the Galerkin equation, we move one factor $\Lambda^s$ onto the data and use~\eqref{kg79}. Its value at the chosen test is then bounded by $C_s D_s E_s(z_j)$. Combining these bounds with coercivity, we obtain
\begin{equation}
E_s(z_j) \leq C_s\bigl(D_s + E_{s-1}(z_j) + (1 + \norm{p}_{s+d_0}) E_{s_0}(z_j)\bigr).
\label{kg82}
\end{equation}

For the finitely many orders through $s_0$, the crude commutator bounds give $E_s(z_j) \leq C_s(D_s + E_{s-1}(z_j))$. An induction therefore yields $E_{s_0}(z_j) \leq CD_{s_0}$. At higher orders we use instead the spectral interpolation inequality
\begin{equation*}
E_{s-1}(z_j) \leq \varepsilon E_s(z_j) + C_{s, \varepsilon} E_0(z_j)
\end{equation*}
which, together with~\eqref{kg82}, proves~\eqref{kg76} with constants uniform in $j$. We can consequently pass to weak limits at each order. By density of the spectral subspaces, the limit is the unique weak solution, and the same estimate proves its smoothness. Notice that the higher coefficient norms appear in the bound, but no smallness is imposed on these norms.

It remains to prove smooth tame dependence. Let $\ZZ^\infty$ denote the space of smooth triples in $\ZZ_0$. Using the formal adjoint of $S$, we define
\begin{equation*}
\mathsf L_p = S^* G^\sharp S + \mathsf R_p: \ZZ^\infty \longrightarrow\ZZ^\infty.
\end{equation*}
We have a smooth tame family of linear maps on a fixed tame space. Coercivity, together with the regularity just proved, shows that each $\mathsf L_p$ is bijective, while estimate~\eqref{kg76} with $H_1 = 0$ gives a tame inverse. To check continuity in the smooth topology, we subtract two weak equations and apply the same estimate to their difference. We can therefore invoke Hamilton's theorem on inverses of linear families~\cite[Part II, Theorem~3.1.1]{Hamilton} to conclude that the inverse depends smoothly tamely on the coefficients. Composition with the source map~$(H_1, H_2) \mapsto S^* H_1 + H_2$ completes the proof.
\end{proof}

Combining Lemma~\ref{kg60}, Proposition~\ref{kg69} and Proposition~\ref{kg75} allows us to solve \eqref{reduced linear eq} for $g\in \HH$. To prove Proposition~\ref{prop right inverse of the linear operator} we still need two more ingredients: to solve the divergence equation and to extend an invariant divergence-free datum $g$ to a function in $\HH$. For this we use two radial homotopies. Given a scalar function~$h$, we define
\begin{equation}
\mathcal Dh = \left(\int_0^1 t\,h(ty, Z)\,dt\right)(y_1 \partial_{y_1} + y_2 \partial_{y_2}).
\label{kg83}
\end{equation}
This gives $\Div \mathcal Dh = h$. For the second construction, let $g = g^1 \partial_{y_1} + g^2 \partial_{y_2} + g^Z U$ be a solenoidal field and set
\begin{equation}
\eta_g = \left[\int_0^1 \bigl(y_2 g^1(ty, Z) - y_1 g^2(ty, Z)\bigr)\,dt\right] dZ+\left[\int_0^1 t\,g^Z(ty, Z)\,dt\right](y_1 dy_2 - y_2 dy_1).
\label{kg84}
\end{equation}
The differential-form homotopy associated with $(y, Z) \mapsto (ty, Z)$ gives $d\eta_g = \iota_g \mu$, because the closed two-form has zero pullback to the central circle. Both homotopies are smooth at the axis and commute with the circle actions. Moreover, differentiating the formulas directly yields
\begin{equation}
\norm{\mathcal Dh}_{C^s_\sigma} \leq C_s \sigma \norm{h}_{C^s_\sigma}, \qquad \norm{\eta_g}_{C^s_\sigma} \leq C_s \sigma \norm{g}_{C^s_\sigma}.
\label{kg85}
\end{equation}
By Sobolev embedding, these bounds give tame estimates on the usual smooth spaces, with a fixed shift in the grading.

Let us extend $\eta_g$ to $M$ with a fixed equivariant extension operator, as in Lemma~\ref{kg26}. Applying $\Pi$ to the metric dual of this extension gives a field~$h_g$, and we define
\begin{equation}
\widetilde g = \CC h_g.
\label{kg86}
\end{equation}
When $\av g = 0$, equivariance ensures that $\av\eta_g = 0$. Since the Hodge projection removes only closed one-forms, the field~$\widetilde g$ coincides with $g$ on the cap and satisfies $\CC^{-1} \widetilde g = h_g$. We have thus obtained a global source for the variational problem without prescribing any boundary condition on $g$.

\begin{proof}[Proof of Proposition~\ref{prop right inverse of the linear operator}]
The maps $\Phi \mapsto G_\Phi, c_\Phi, \av c_\Phi$ are smooth tame. We extract the profiles of $\av c_\Phi$ using Lemma~\ref{lem:profiles} and extend the coefficients by the construction above. After shrinking one finite-order neighborhood of $\Phi_0$, the resulting coefficients belong to the neighborhood required in Proposition~\ref{kg75}.

Given $(h, g) \in \mathcal Y^\infty$, we use $\mathcal Dh$ as a particular solution of the divergence equation and put
\begin{equation*}
g_1 = g - (M_\Phi \mathcal Dh)_2.
\end{equation*}
This second datum is solenoidal and jointly mean-free, since the second component of $M_\Phi$ takes values in that space. We apply~\eqref{kg86} to $g_1$ and solve~\eqref{kg74} with data~$(h_{g_1}, 0)$. Proposition~\ref{kg75} supplies a smooth solution $z$. Applying Lemma~\ref{kg60} to $\widetilde g_1 = \CC h_{g_1}$, we then obtain the solenoidal correction to the particular solution. The full right inverse is therefore
\begin{equation}
R_\Phi(h, g) = \mathcal Dh + \left.(T_b e + T_j \mathsf P\mathfrak r)\right|_{\Omega_{\sigma, N}}, \qquad \mathfrak r = (\Id - \av) \curl_\mu(G(Sz, \cdot)) - \CC h_{g_1}.
\label{kg87}
\end{equation}
It satisfies $M_\Phi R_\Phi(h, g) = (h, g)$.

Each operation in this formula has a fixed loss of derivatives and depends smoothly tamely on $\Phi$. For completeness, the denominator occurring in the reconstruction by current transport can be written as
\begin{equation*}
\lambda^{-1} = \frac1m \left(\beta + \alpha\omega\frac{k}{m}\right)^{-1}.
\end{equation*}
On $\mathsf P$, we have $\abs{m}^{-1} \leq N$ and $\abs{k/m} \leq N\delta$, and the expression in parentheses is uniformly separated from zero. Its coefficient derivatives therefore satisfy uniform product estimates as the modal ratio $k/m$ ranges over a fixed bounded interval. These coefficients preserve the weights, so orthogonal summation yields the Sobolev bounds for $T_j$. This proves the asserted right-inverse property together with the tame estimates.
\end{proof}

For clarity, the spaces in this construction have the smoothing operators required in the tame category. On $M$, spectral smoothing commutes with the Hodge and angular projections. On the physical cap, the usual smooth scale is obtained by extension and restriction, and its solenoidal, jointly mean-free target is the range of the tame projection
\begin{equation*}
v \longmapsto (\Id - \av)(v - \mathcal D\Div v).
\end{equation*}
Consequently, all the spaces in the linear construction are fixed tame spaces, including at the magnetic axis.

We finish the section by proving the estimate that will allow us to move the axis. At the reference immersion, the coefficients are all $U$-invariant, and therefore every longitudinal sector is preserved by the construction. On a finite collection of sectors we can choose the band so as to isolate the exact resonances. The augmented reconstruction is then of order zero, and this yields the factor of the thickness needed below.

\begin{corollary}\label{kg88}
Let $\mathcal M \subset N^{-1} \Z\setminus\{0\}$ be finite and satisfy $\kappa_* m \neq \pm1$ for $m \in \mathcal M$. Choose the band so that
\begin{equation*}
\delta < \tfrac12\min_{m\in\mathcal M} \inf\{\abs{n - \kappa_* m}: n \in \Z, \ n - \kappa_* m \neq 0\}.
\end{equation*}
For a smooth solenoidal physical datum~$g$ whose longitudinal frequencies lie in $\mathcal M$, $\xi_{\rm corr} = R_{\Phi_0}(0, g)$ is solenoidal, has the same sectors, and satisfies
\begin{equation}
\begin{gathered}
\curl[b, \xi_{\rm corr}] - [j_0, \xi_{\rm corr}] = g, \\
\norm{[b, \xi_{\rm corr}]}_{C^r_\sigma} + \norm{\xi_{\rm corr}^y(0, \cdot)}_{C^r(\T_N)} \leq C_r \sigma \norm{g}_{C^{r+r_*}_\sigma},
\end{gathered}
\label{kg89}
\end{equation}
where $\xi_{\rm corr}^y$ denotes the two transverse Cartesian components and $r_*$ is fixed. The constants are uniform for sufficiently small~$\sigma$ and may depend on the fixed set~$\mathcal M$.
\end{corollary}

\begin{proof}
For each of the sectors under consideration, the selected weights satisfy $k = 0$. Hence $J_1 = 0$ and $Sz = e + J_2 \tau$, so the augmented form is an order-zero system in $(e, f, \tau)$. By Proposition~\ref{kg69}, its positive lower bound is uniform in $\sigma$. We emphasize that the scalar unknown $f$ and its positive quadratic term are still present at these exact resonances.

We apply the energy argument of Proposition~\ref{kg75}, now using the normalized sphere grading $\Lambda_{R_\Sigma} = (1 + R_\Sigma^2 \Delta_\Sigma)^{1/2}$ and the norm
\begin{equation*}
\norm{u}_{s, R_\Sigma} = R_\Sigma^{-1} \norm{\Lambda_{R_\Sigma}^s u}_{L^2(M)}.
\end{equation*}
Here $\Delta_\Sigma$ is the partial Hodge Laplacian on the relevant scalar or one-form bundle. The gradings commute with the Hodge and angular projections, as well as with $Q_E$ and $\CC^{-1}$. For a fixed longitudinal frequency, they define the Sobolev norms of the unit sphere in physical tensor components. By Lemma~\ref{kg26}, the reference coefficients have uniform rescaled bounds, including $\omega'$ and $\beta_0'$, and the same holds for $d\varrho/R_\Sigma$ and $X/R_\Sigma$. Furthermore,~\eqref{kg50} gives uniform bounds for $R_\Sigma^{-1} \CC^{-1}$ and $R_\Sigma^{-1} \CC^{-1} U$ at fixed~$m$. To use these bounds in~\eqref{kg51}, we move $U$ past the invariant coefficients and place it next to $\CC^{-1}$. The positive gap between the nonzero values of $k$ and zero then bounds $A_E$ and its commutators with $\Lambda_{R_\Sigma}^s$ uniformly. On the selected weights, the operators $J_2$, $\mathsf B$ and $\mathsf C$ have order zero with uniformly bounded rescaled coefficients. We can therefore commute the augmented equation and interpolate the lower-order terms, obtaining
\begin{equation*}
\norm{e}_{s, R_\Sigma} + \norm{f}_{s, R_\Sigma} + \norm{\tau}_{s, R_\Sigma} \leq C_s \norm{h_g}_{s, R_\Sigma}, \qquad \norm{w}_{s, R_\Sigma} \leq C_s \norm{h_g}_{s, R_\Sigma}, \quad w = Sz.
\end{equation*}
Since the Hodge projection is a contraction and commutes with the grading,~\eqref{kg85} and the rescaled extension bounds imply $\norm{h_g}_{s, R_\Sigma} \leq C_s \sigma \norm{g}_{C^s_\sigma}$. We now apply Sobolev embedding on the unit sphere, with a fixed shift $r_* > 2$, to obtain the stated estimate for $w = [b, \xi_{\rm corr}]$. Every operation used in this argument preserves the longitudinal sectors.

It only remains to control the normal displacement of the axis. In a sector $e^{\ii mZ}$, the transverse bracket on the axis is given exactly by
\begin{equation*}
[b, \xi_{\rm corr}]^y(0, Z) = \omega(0)(-\ii m\kappa_* \Id - J) \xi_{\rm corr}^y(0, Z), \qquad \det(-\ii m\kappa_* \Id - J) = 1 - m^2 \kappa_*^2.
\end{equation*}
By our hypothesis, all these finitely many matrices are invertible. The bound already obtained for $w$ therefore controls the normal displacement as well and proves the remaining term in~\eqref{kg89}.
\end{proof}

\section{Braiding the magnetic axis}\label{kg90}

Our objective is to find a smooth curve of exact solutions for which the first variation of the magnetic axis is a braid of the prescribed type. There are two conditions to meet: the variation must belong to the kernel of the linearized equation, and its normal trace must have the required knot type. For this, we choose the braid before fixing the thickness, so that its finite set of longitudinal sectors avoids the two resonances on the axis. Corollary~\ref{kg88} will then allow us to correct a trial displacement by a term of size $O(\sigma)$ on the axis, without losing the braid. Once this kernel vector has been constructed, we use Hamilton's implicit-function theorem to integrate it, taking care to preserve its full first variation. This gives the family in~\eqref{kg105}; in the next section we will restrict these immersions to domains on which they are embeddings.

Recall that Alexander's theorem~\cite{Alexander} represents any knot as the connected closure of a braid around a line. For a closed braid with $N$ strands, this amounts to choosing a smooth profile in the normal plane,
\begin{equation*}
\mathbf F(\varphi) = (f(\varphi), g(\varphi)), \qquad \varphi \in \T_N,
\end{equation*}
such that, for each $\theta \in \T$, the points
\begin{equation}
\mathbf F(\theta + 2\pi j), \qquad 0 \leq j < N,
\label{kg91}
\end{equation}
are pairwise distinct. By compactness, the distance between any two strands has a positive minimum, and a transverse compression places the braid inside an arbitrarily thin tube around a circle. We can also take $N$ to be any larger integer, without affecting the knot type: the addition of one strand with one crossing produces, after taking the closure, only a local detour which can be removed by an isotopy.

It is convenient to subtract the motion common to all the strands in a fibre. To this end, we introduce the average
\begin{equation*}
\deckav \mathbf F(\varphi) = \frac1N\sum_{j=0}^{N-1} \mathbf F(\varphi + 2\pi j).
\end{equation*}
The replacement of $\mathbf F$ by $(\Id - \deckav) \mathbf F$ does not change the differences between any two points in~\eqref{kg91}. In particular, the homotopy obtained by translating the strands preserves the braid, and, after a transverse compression, it remains in the half-plane where the cylindrical radius is positive. We shall call a profile with this normalization \emph{centered}. The longitudinal Fourier frequencies of a centered profile satisfy
\begin{equation}
m = j/N, \qquad j \in \Z\setminus N\Z.
\label{kg92}
\end{equation}
Thus the Fourier expansion contains neither the zero frequency nor any of the integer frequencies.

\begin{lemma}\label{kg93}
For every knot type and every $\kappa > 0$, there are $N \geq 2$ with $\kappa_* := N\kappa > 1$ and a centered real trigonometric polynomial $\mathbf F$ representing that knot such that
\begin{equation}
\kappa_* m \neq \pm1
\label{kg94}
\end{equation}
for every longitudinal frequency~$m$ of $\mathbf F$.
\end{lemma}
\begin{proof}
Start with a braid of the prescribed knot type and stabilize it until $N\kappa > 1$. After centering the profile, we approximate it by a centered trigonometric polynomial. In view of $m = j/N$, the frequencies to be excluded are characterized by $\kappa j = \pm1$, so there is nothing to remove unless $1/\kappa$ is a positive integer. We are therefore left with the following elementary question: how can one remove a finite set of longitudinal Fourier modes while preserving the braid?

For this purpose we identify the normal plane with $\C$, so that the profile is written as $f + \ii g$, and consider the smooth periodic function
\begin{equation*}
\mathfrak u(\varphi) = \left(\frac{1 - 3^{-1/2} e^{-\ii\varphi}}{1 - 3^{-1/2} e^{\ii\varphi}}\right)^2
\end{equation*}
which has modulus one, winding number zero and zero mean. To check the last assertion, expand the denominator in a geometric series: its constant Fourier coefficient is $1 - 4/3 + 3/9 = 0$. Moreover, replacing $3^{-1/2}$ by $t3^{-1/2}$, with $0 \leq t \leq 1$, joins the constant function~$1$ to $\mathfrak u$ through periodic functions of modulus one, and also proves that the winding number is zero. Now choose a positive integer~$\ell_0$ large enough that translation by any nonzero multiple of $\ell_0$ moves the finite Fourier support of $f + \ii g$ away from $m = \pm\kappa_*^{-1}$. Multiplying the profile by $\mathfrak u(\ell_0 \varphi)$ removes the modes that we wish to exclude, since the multiplier has zero constant coefficient and all its other frequencies are nonzero multiples of $\ell_0$.

Geometrically, this multiplication rotates every strand in a given angular fibre by the same angle. Hence the homotopy just described preserves both the separation of the strands and the knot type. Notice also that the multiplier is $2\pi$-periodic, so the centered condition is preserved. A sufficiently close Fourier truncation of the rotated profile now gives a trigonometric polynomial with the same braid type, and the coefficients at the excluded frequencies are still zero.
\end{proof}

For the asymmetry assertion we make one further choice in the profile. Since $-(f'' + f)$ is the first variation of the curvature of a circle, up to a positive constant, a unique nondegenerate global maximum of this function will single out a point of the deformed axis. Its jets at that point will then distinguish the two possible directions along the curve.

\begin{lemma}\label{kg95}
The profile in Lemma~\ref{kg93} can be chosen so that
\begin{equation*}
K_{\rm lin} = -(f'' + f)
\end{equation*}
has a unique nondegenerate global maximum at $\varphi_*$ and
\begin{equation}
K_{\rm lin}'''(\varphi_*) \neq 0, \qquad g'''(\varphi_*) + g'(\varphi_*) \neq 0.
\label{kg96}
\end{equation}
This can be achieved with an arbitrarily small~$C^1$ change of the initial centered braid.
\end{lemma}
\begin{proof}
We first introduce the Fourier projection~$\mathsf Q$ onto those longitudinal modes~$m = \pm\kappa_*^{-1}$ which belong to $N^{-1} \Z$, with the convention $\mathsf Q = 0$ when neither mode occurs. As $\kappa_* > 1$, these are noninteger frequencies and $\mathsf Q\deckav = 0$. On centered functions on $\T_N$, the curvature operator~$L_{\rm curv} = -(\partial_\varphi^2 + 1)$ is invertible: indeed, the only frequencies that could obstruct its inversion are $\pm1$, and these have already been removed. The periodic Green kernel with these two frequencies removed is continuous, and its one-sided first derivatives are bounded. This follows by solving the smooth second-order ordinary differential equation away from the origin and observing that the first derivative has only a finite jump there. Convolution with this kernel therefore yields
\begin{equation}
\norm{L_{\rm curv}^{-1} H}_{C^1} \leq C_N \norm{H}_{L^1} \qquad (\deckav H = 0).
\label{kg97}
\end{equation}
Since this inverse commutes with $\deckav$ and $\mathsf Q$, applying it cannot destroy either the centered condition or~\eqref{kg94}.

Choose a nonnegative smooth function~$\chi$ with support in $(-1, 1)$, whose unique maximum is at zero, and such that $\chi'(0) = 0$, $\chi''(0) < 0$ and $\chi'''(0) \neq 0$. Fix a point~$\varphi_0$. For $\varepsilon > 0$ small, we periodize the resulting narrow bump and define
\begin{equation*}
H_\varepsilon = \varepsilon^{-1/2}(\Id - \deckav - \mathsf Q) \chi\bigl((\varphi - \varphi_0)/\varepsilon\bigr).
\end{equation*}
For sufficiently small support, the deck translates of the bump are disjoint. Its projection under $\mathsf Q$, on the other hand, is $O(\varepsilon^{1/2})$ in every fixed smooth norm. Thus, apart from this small correction, $H_\varepsilon$ consists of a positive bump at $\varphi_0$, with coefficient~$1 - 1/N$, and negative bumps at the remaining $N - 1$ points in its deck orbit. The inverse estimate~\eqref{kg97} gives
\begin{equation*}
\norm{L_{\rm curv}^{-1} H_\varepsilon}_{C^1} \leq C_N \varepsilon^{1/2}.
\end{equation*}
Replacing $f$ by $f + L_{\rm curv}^{-1} H_\varepsilon$ changes the curvature variation to $K_{\rm lin} + H_\varepsilon$. The maximum of this new function must lie in the positive bump, whose height tends to infinity, whereas the perturbation outside it is bounded above by $O(\varepsilon^{1/2})$. Near the maximum of $\chi$, the negative second derivative of the bump dominates the fixed background. Outside a smaller neighborhood of that maximum, there is instead a fixed relative gap in the height of the bump. These two observations show that the global maximum is unique and nondegenerate. Its position differs from $\varphi_0$ by $O(\varepsilon^{5/2})$, so the leading term of the third derivative at the maximum is
\begin{equation*}
(1 - 1/N) \varepsilon^{-7/2} \chi'''(0) \neq 0.
\end{equation*}
This proves the first condition in~\eqref{kg96}.

We obtain the other required jet by adding a sufficiently localized bump to $g$ near this maximum and projecting it by $\Id - \deckav - \mathsf Q$. Such a bump can prescribe a nonzero change of $g''' + g'$ at the chosen point. The subtraction of finitely many Fourier modes does not cancel this change: its derivatives are bounded in terms of the $L^1$ norm of the bump, while the localized third derivative has the larger scale. We may choose the amplitude arbitrarily small, avoiding the single value that would cancel the jet already present. In particular, the changes of $f, g$ are small in $C^1$, and the strands remain separated.

It remains to return to trigonometric polynomial profiles. We approximate the modified centered profile by centered real trigonometric polynomials in $C^8$, keeping the excluded Fourier coefficients equal to zero. If the approximation is sufficiently close, the unique nondegenerate maximum, the two nonzero jets and the separation of the strands all persist. Consequently~\eqref{kg94} is preserved as well.
\end{proof}

The profile now has the prescribed knot type and the local geometric properties that will rule out symmetries. Our next task is to realize it, up to a small error, as the normal trace of a field in the kernel. Let $\mathcal M$ denote the finite symmetric set of longitudinal weights appearing in $\mathbf F(-Z)$. We first fix $L$ sufficiently large for Proposition~\ref{kg69}, and then choose $\delta$ small enough for that proposition and for the additional requirement
\begin{equation}
\delta < \frac12\min_{m\in\mathcal M} \inf\{\abs{n - \kappa_* m}: n \in \Z, \ n - \kappa_* m \neq 0\}.
\label{kg98}
\end{equation}
Indeed, in a sector $\mathcal L_U = \ii m$ the magnetic divisor is $k = n - \kappa_* m$, so this choice ensures that the selected projection~$\mathsf P$ contains only the exact resonances in the sectors under consideration. In particular, $J_1 = 0$ in these sectors, while~\eqref{kg94} guarantees that the transverse bracket on the axis is invertible. The full right inverse continues to treat every other mode as before; only its use on this finite collection of sectors has been specified more precisely. We can now choose the thickness~$\sigma$ small enough for both the inverse estimates and the geometric requirements below. After making these choices, we fix the coefficient neighborhood and the deformation interval.

Using the volume coordinates of Section~\ref{sec reference}, consider the exactly solenoidal trial field
\begin{equation}
\xi_{\rm tr} = f(-Z) \partial_{y_1} + g(-Z) \partial_{y_2}.
\label{kg99}
\end{equation}
To apply the small-axis estimate, we need its defect in the linearized equation to remain bounded as the thickness tends to zero. This follows by examining a single Fourier component: since $b = \omega V$, we have
\begin{equation}
[b, e^{\ii mZ} \partial_{y_j}] = e^{\ii mZ} \bigl[\omega(-\ii m\kappa_* \Id - J) e_j - \omega'y_j V\bigr].
\label{kg100}
\end{equation}
The leading transverse coefficient is independent of $y$ up to $O(\sigma^2)$, and the remaining terms are $O(\sigma)$ in the rescaled component norms. The metric has the same property, with a constant leading shear and an $O(\sigma)$ remainder by~\eqref{kg24}, so one physical transverse derivative of either remainder is bounded. Longitudinal differentiation introduces only the fixed frequencies in $\mathcal M$. Moreover, the current bracket is
\begin{equation*}
[j_0, e^{\ii mZ} \partial_{y_j}] = \ii m\beta_0 e^{\ii mZ} \partial_{y_j} - \beta_0'y_j e^{\ii mZ} U.
\end{equation*}
and therefore the defect $\curl[b, \xi_{\rm tr}] - [j_0, \xi_{\rm tr}]$ is uniformly bounded in every fixed rescaled norm. Since it is solenoidal and its longitudinal sectors lie in $\mathcal M$, it has the properties required for the correction.

Corollary~\ref{kg88} now applies in each of these sectors. Summing the corrections gives a field~$\xi_{\rm corr}$ which is real, because all the operations commute with complex conjugation. Thus, setting $\xi = \xi_{\rm tr} - \xi_{\rm corr}$, we obtain
\begin{equation}
M_{\Phi_0} \xi = 0, \qquad \xi_{\rm corr}^y(0, \cdot) = O(\sigma)
\label{kg101}
\end{equation}
in each fixed smooth norm along the axis. Recall here that the mean-zero angular primitive $u$ vanishes on the axis, so that $Z = -\varphi$ there. The normal components of $D\Phi_0 \partial_{y_1}$ and $D\Phi_0 \partial_{y_2}$ form the matrix $\Id + O(\sigma)$, whereas $D\Phi_0 U$ is tangential; consequently the normal trace of $D\Phi_0 \xi$ is a centered profile with finite Fourier support and differs from $\mathbf F$ by $O(\sigma)$ in $C^s(\T_N)$ for any fixed~$s$. It follows that, on choosing $\sigma$ small enough, both strand separation and the conditions of Lemma~\ref{kg95} survive the correction. This choice involves only finitely many norms. Henceforth $(f, g)$ will denote the corrected normal trace, $K_{\rm lin} = -(f'' + f)$ its curvature variation, and $\varphi_*$ the point of maximum of that variation.

It remains to integrate the kernel vector into exact immersions, preserving the entire first variation that we have prescribed. In the linearization identity~\eqref{kg37}, take $\xi = R_\Phi h$. The right-inverse relation $M_\Phi R_\Phi h = h$, supplied by Proposition~\ref{prop right inverse of the linear operator}, then gives
\begin{equation}
D\mathcal F(\Phi)[D\Phi R_\Phi h] = h + \err(\mathcal F(\Phi), R_\Phi h).
\label{kg102}
\end{equation}
The maps $D\Phi R_\Phi h$ and $\err(r, R_\Phi h)$ are smooth tame in all their arguments. Moreover, the error is bilinear in $(r, h)$ when $\Phi$ is fixed, which is precisely the quadratic-error condition in Hamilton's implicit-function theorem \cite[Part III, Theorem 3.3.1]{Hamilton}.

We shall use the following consequence of that theorem to retain the prescribed tangent vector.

\begin{lemma}\label{kg103}
Let $\mathcal X, \mathcal Y$ be tame Fr\'echet spaces, $\mathcal O \subset \mathcal X$ an open set, and $F: \mathcal O \to \mathcal Y$ a smooth tame map with $F(x_0) = 0$. Suppose that, on a neighborhood of $x_0$, there are smooth tame maps $Q_x: \mathcal Y \to \mathcal X$, linear in their argument, and $B_x: \mathcal Y \times \mathcal Y \to \mathcal Y$, bilinear in their arguments, such that
\begin{equation*}
DF(x) Q_x h = h + B_x(F(x), h).
\end{equation*}
For every fixed~$v \in \ker DF(x_0)$ there is a smooth curve of zeros
\begin{equation}
x(t) = x_0 + tv + O(t^2)
\label{kg104}
\end{equation}
in every seminorm.
\end{lemma}
\begin{proof}
Set~$x = x_0 + t(v + w)$ and define, for $t \neq 0$,
\begin{equation*}
A(t, w) = \frac{F(x_0 + t(v + w))}{t}.
\end{equation*}
Using the identity
\begin{equation*}
A(t, w) = \int_0^1 DF(x_0 + ut(v + w))(v + w)\,du
\end{equation*}
we see that $A$ extends smoothly and tamely across $t = 0$, where $A(0, w) = DF(x_0) w$ and thus $A(0, 0) = 0$. Differentiation in the unknown gives
\begin{equation*}
D_w A(t, w) = DF(x_0 + t(v + w)), \qquad D_w A(t, w) Q_x h = h + tB_x(A(t, w), h).
\end{equation*}
Hamilton's theorem with unknown $w$ and parameter~$t$ provides a smooth solution $w(t)$. To retain the prescribed tangent, we use the iteration in the proof of the theorem \cite[p.~212]{Hamilton} with initial value $w_0(t) = 0$. Every correction contains the residual $A(t, w_j(t))$, and since $A(0, 0) = 0$, induction gives $w_j(0) = 0$ for every $j$. The limiting solution therefore satisfies $w(0) = 0$, so its smoothness yields $w(t) = O(t)$. Substitution in the ansatz proves~\eqref{kg104}.
\end{proof}

At the reference immersion, the identities~\eqref{kg37} and~\eqref{kg101} imply $D\mathcal F(\Phi_0)[D\Phi_0 \xi] = 0$. Thus the kernel direction just constructed satisfies the hypothesis of the lemma. We use the source and target spaces, together with the single coefficient neighborhood, furnished by Proposition~\ref{prop right inverse of the linear operator}. The lemma then gives a jointly smooth family of orientation-preserving immersions of the fixed covered source, satisfying
\begin{equation}
\mathcal F(\Phi_t) = 0, \qquad \Phi_t = \Phi_0 + tD\Phi_0 \xi + O(t^2).
\label{kg105}
\end{equation}
The rank remains three because this is an open condition in $C^1$. By construction, the normal component of the first variation is exactly the braid profile~$(f, g)$ chosen above.

\section{Proof of the main theorem}\label{kg106}

To prove Theorem~\ref{kg3} and Corollary~\ref{kg6}, we use the exact immersions constructed in Section~\ref{kg90}. Their first variation already contains the prescribed braid; the point is to show that the same knot appears on the axis of an exact solution and that a sufficiently small surrounding tube is embedded. On this tube the pressure and magnetic dynamics are supplied by the construction. The additional assertions will follow from the geometric conditions imposed on the axis.

Write the reference axis as $R_\sigma e_R(\varphi) + z_\sigma e_z$, with $R_\sigma > 0$. The expansion in~\eqref{kg105} gives its deformed image as
\begin{equation}
\Gamma_t(\varphi) = R_\sigma e_R(\varphi) + z_\sigma e_z + t\bigl(f(\varphi) e_R(\varphi) + h(\varphi) e_\varphi(\varphi) + g(\varphi) e_z\bigr) + O(t^2).
\label{kg107}
\end{equation}
Here the pair~$(f, g)$ is the normal trace fixed in Section~\ref{kg90}. Observe that the physical angular coordinate has a lift of the form~$\varphi + t\,h(\varphi)/R_\sigma + O(t^2)$ with positive derivative and period increment $2\pi N$. We can therefore use it to reparametrize the curve and divide the radial and vertical displacements from the reference circle by $t$, obtaining a profile that converges to $(f, g)$ in $C^1$ as $t\downarrow0$. As the limiting strands have a positive minimum separation, straight interpolation to this limiting profile stays among embedded closed braids for small positive $t$. Hence $\Gamma_t$ is an embedded knot of the required type.

The immersion property then implies that $\Phi_t$ is injective on a sufficiently small neighborhood of its axis. Indeed, if no such neighborhood existed, we could find distinct source points~$x_j, y_j$ approaching the axis with $\Phi_t(x_j) = \Phi_t(y_j)$. After passing to subsequences, their limits would be points~$x, y$ on the axis, and the injectivity there would give $x = y$, contradicting local invertibility for large~$j$. In particular, we can choose $0 < \varrho_t < \varrho_b$ such that $\Phi_t$ embeds $\{\varrho \leq \varrho_t\}$, and set
\begin{equation*}
\Omega_t = \Phi_t(\{\varrho < \varrho_t\}), \qquad B_t \circ \Phi_t = D\Phi_t\,b.
\end{equation*}
For each fixed positive~$t_*$, the same closed source sublevel works for nearby parameters, by openness of embeddings on a compact source. The final radius need not be bounded below as $t\downarrow0$; it is allowed to tend to zero.

By Proposition~\ref{kg33} applied to the immersion~$\Phi_t$ with $a = \varrho_t$, the exact equation~\eqref{kg105} and the injectivity just proved allow us to use~\eqref{kg34} to define a single pressure $P_t$ on $\overline\Omega_t$. Thus $(B_t, P_t)$ satisfies~\eqref{kg2}. The same proposition gives positivity in the interior and a nondegenerate maximum precisely along $\Gamma_t$, with regular pressure levels forming the required foliation of the complement by tori.

On each regular pressure torus put $\zeta = -Z/N \in \T$. Since the field is transported by the immersion, taking $s$ to be physical time gives the flow
\begin{equation}
(\Theta, \zeta) \longmapsto (\Theta, \zeta) + s\,\omega(\varrho)(1, \kappa) \pmod{2\pi}.
\label{kg108}
\end{equation}
where~$\kappa_*/N = \kappa$ by the choice in Lemma~\ref{kg93}. If $\kappa$ is irrational, every nonconstant Fourier character has nonzero frequency and hence a time average tending uniformly to zero. Approximation by trigonometric polynomials gives the same uniform convergence to the torus average for every continuous function, proving unique ergodicity. The invariant probability measure has full support, so every orbit is dense. If $\kappa = A/d$ is in lowest terms, the regular trajectories instead have winding class $(d, A)$ and minimal period~$2\pi d/\omega(\varrho)$, and the axis has minimal period $2\pi/(\kappa\omega(0))$. These are the periods in physical time, retaining their dependence on the pressure level.

Recall that the longitude transported by the immersion can differ from the preferred longitude of the knot. With consistent orientations, it is homologous to the preferred zero-linking longitude plus $\ell$ times the transported meridian, for an integer~$\ell \in \Z$. In the preferred marking the signed longitudinal-to-meridional slope is consequently
\begin{equation}
\frac{\kappa}{1 + \ell\kappa},
\label{kg109}
\end{equation}
with the projective interpretation if the denominator vanishes. Since the integer~$\ell$ is constant along an embedded branch, this also gives the change of marking described after Theorem~\ref{kg3} and completes its proof.

For Corollary~\ref{kg6}, the additional choice of profile in Lemma~\ref{kg95} will be used through the following observation about curves.

\begin{lemma}\label{kg110}
Let a smooth closed embedded curve have positive curvature. Suppose that its curvature has a unique global maximum at $x_*$, the third arclength derivative of curvature there is nonzero, and its torsion at $x_*$ is nonzero. Then its Euclidean isometry group is trivial.
\end{lemma}
\begin{proof}
Any isometry preserving the curve fixes $x_*$, the unique point of maximum curvature, so the restriction of the isometry in arclength coordinates based at $x_*$ is either $s \mapsto s$ or $s \mapsto -s$. In the second case the curvature would be even, which is impossible because its third derivative is nonzero. Thus the curve is fixed pointwise, and the orthogonal linear part fixes the first three derivatives of an arclength parametrization at $x_*$. These span $\R^3$ by positive curvature and nonzero torsion, so the linear part is the identity. The translation part must also vanish, since the isometry has a fixed point.
\end{proof}

For the curvature $K_t$ and torsion $\tau_t$ of the curve in~\eqref{kg107}, substitution in the usual formulas gives the first variations
\begin{equation}
\dot K_0 = -\frac{f'' + f}{R_\sigma^2} = \frac{K_{\rm lin}}{R_\sigma^2}, \qquad \dot\tau_0 = \frac{g''' + g'}{R_\sigma^2}.
\label{kg111}
\end{equation}
Indeed, the tangential and vertical contributions cancel in $K = \abs{\Gamma' \times \Gamma''}/\abs{\Gamma'}^3$, while in the torsion formula the numerator $\det(\Gamma', \Gamma'', \Gamma''')$ is $tR_\sigma^2(g''' + g') + O(t^2)$ and the denominator is $R_\sigma^4 + O(t)$. For small~$t$ the curvature therefore stays positive. Its maximum is unique and occurs at $\varphi_t = \varphi_* + O(t)$, because the rescaled variation $(K_t - R_\sigma^{-1})/t$ converges smoothly to $K_{\rm lin}/R_\sigma^2$, whose maximum was chosen to be unique and nondegenerate. At this point,
\begin{equation}
\frac{d^3 K_t}{ds_t^3} = \frac{t}{R_\sigma^5} K_{\rm lin}'''(\varphi_*) + O(t^2), \qquad \tau_t = \frac{t}{R_\sigma^2}(g''' + g')(\varphi_*) + O(t^2).
\label{kg112}
\end{equation}
and both leading coefficients are nonzero by construction. Lemma~\ref{kg110} therefore applies for all sufficiently small positive $t$.

Assume that an isometry preserves the equilibrium, possibly reversing $B_t$. The force equation makes the transformed pressure differ from the original one by an additive constant; this constant is zero because the pressure vanishes on the boundary. The same argument applies to an orientation-reversing isometry, as the two orientation signs from curl and the cross product cancel. Such an isometry must therefore preserve the unique critical curve of the pressure, and its trivial isometry group forces the ambient isometry to be the identity.

To separate the members of the family up to similarities, let $\ell(C)$ be the length of a curve $C$ and consider
\begin{equation}
\mathcal W(C) = \ell(C) \int_C(K - \overline K)^2\,ds, \qquad \overline K = \ell(C)^{-1} \int_C K\,ds.
\label{kg113}
\end{equation}
This expression is invariant under spatial similarities and reparametrizations. As a parameter integral it is smooth in $t$ through the multiply covered circle $\Gamma_0$, where the curvature variance vanishes. Since the curvature first variation is the nonconstant function~$K_{\rm lin}/R_\sigma^2$, expansion at this circle yields
\begin{equation}
\mathcal W(\Gamma_t) = ct^2 + O(t^3), \qquad c > 0.
\label{kg114}
\end{equation}
The dependence on $t$ being smooth, we also have $\frac d{dt} \mathcal W(\Gamma_t) = 2ct + O(t^2) > 0$ on a sufficiently small positive interval, so axes with distinct parameters in that interval cannot be similar.

An equivalence by a spatial similarity, a nonzero constant field rescaling and an additive pressure constant sends the unique critical axis of one equilibrium to that of the other. It must therefore preserve~\eqref{kg113}, which contradicts the strict variation at distinct parameters. Near each fixed positive parameter the source sublevel can be kept unchanged, as noted above, giving a smooth local family of inequivalent equilibria. This proves Corollary~\ref{kg6}.

\section{Equilibria with nonconstant rotational transform}\label{kg115}

We finish by adapting the construction to a small nonconstant rotational transform, without changing the prescribed knot type. The reference equilibrium in this argument is still smooth. The restriction to finite regularity enters when deforming it in three dimensions: differentiating the variational equation now produces terms of the same order as the unknown. Their coefficients tend to zero with the twist, however, so the terms can be absorbed at any finite collection of derivative orders fixed beforehand. This is the point at which the size of the allowed twist will depend on the desired regularity.

\begin{proposition}\label{kg116}
Fix an integer~$r \geq 2$, a number~$\kappa > 0$, and a knot type. There are a bounded~$C^{r+1}$ solid torus~$\Omega$ of the prescribed type and fields~$B, P$ of class $C^r$ up to the boundary satisfying~\eqref{kg2}, with $B$ nowhere zero on $\overline\Omega$ and the pressure properties of Theorem~\ref{kg3}\textnormal{(i)}. On every regular pressure torus the magnetic flow is $C^r$-conjugate, with time unchanged, to
\begin{equation*}
\omega_\varepsilon(\varrho) \left(\partial_\Theta + \kappa_\varepsilon(\varrho) \partial_\zeta\right), \qquad \omega_\varepsilon > 0, \quad \kappa_\varepsilon > 0, \quad \partial_\varrho \kappa_\varepsilon < 0,
\end{equation*}
for a fixed~$\varepsilon > 0$. The smooth reference profiles approach their constant-slope counterparts as $\varepsilon\downarrow0$ at fixed thickness, with $\kappa_\varepsilon \to \kappa$. The choice of $\varepsilon$ may depend on $r$, $\kappa$, and the braid.
\end{proposition}

\begin{proof}
Choose the covering degree and the braid as in Lemma~\ref{kg93}, so that $\kappa_* = N\kappa > 1$. We begin with a perturbation of the disk solution from Lemma~\ref{kg13}. Let a fixed continuous linear smooth extension operator act on $r_F = (F_\sigma - 2)/\sigma^2$, with its value at $\sigma = 0$ defined by smooth continuation, and define the extended profile to be $2 + \sigma^2\operatorname{Ext}r_F$. This procedure preserves joint smoothness in the parameters and the expansion $F_\sigma = 2 + O(\sigma^2)$ on a fixed neighborhood of $[0, 1/2]$. Using again the notation $F_\sigma$ for this extension, we consider the Dirichlet problem
\begin{equation}
L_\sigma \Psi_{\sigma, \varepsilon} = F_\sigma(\Psi_{\sigma, \varepsilon}) + 8\varepsilon\Psi_{\sigma, \varepsilon}, \qquad \Psi_{\sigma, \varepsilon} |_{\partial D} = \tfrac12.
\label{kg117}
\end{equation}
At $(\sigma, \varepsilon) = (0, 0)$, the Dirichlet linearization is $\Delta$. Thus the implicit function theorem and elliptic regularity provide a smooth strictly convex solution, with parameter estimates uniform through any prescribed finite order. For $\varepsilon = 0$, local uniqueness identifies it with the solution $\Psi_\sigma$ of Lemma~\ref{kg13}. In planar area coordinates it can be written as $\Psi_{\sigma, \varepsilon} = H_{\sigma, \varepsilon}(\mathcal I)$, with $\mathcal I$ equal to enclosed area divided by $2\pi$. Let us check that these coordinates are smooth at the minimum. In Morse coordinates, denote half the squared radius by $\mathcal I_0$ and write the area form as $J_{\rm area}\,d\mathcal I_0 \wedge d\theta$. The changes $\mathcal I = \int_0^{\mathcal I_0} \langle J_{\rm area}\rangle_\theta(s)\,ds$ and $\theta_1 = \theta + \partial_\theta^{-1}(J_{\rm area}/\langle J_{\rm area}\rangle_\theta - 1)$ give $d\mathcal I \wedge d\theta_1 = J_{\rm area}\,d\mathcal I_0 \wedge d\theta$. Since the ratio $\mathcal I/\mathcal I_0$ and the mean-zero angular primitive are smooth in Cartesian coordinates, this construction extends smoothly through the centre. We then continue along the regular level curves, choosing the angular origin on each level to agree with Section~\ref{sec reference} when $\varepsilon = 0$. In particular, at zero twist we have $H_{\sigma, 0}(\mathcal I) = \mathcal I$. If $\sigma = 0$, uniqueness makes the solution on the disk radial, so $\mathcal I = (x^2 + y^2)/2$ and $\Delta(\mathcal I^2) = 8\mathcal I$. Differentiation of~\eqref{kg117} consequently yields
\begin{equation*}
H_{0, \varepsilon}(\mathcal I) = \mathcal I + \varepsilon(\mathcal I^2 - \tfrac14) + O(\varepsilon^2), \qquad \partial_{\mathcal I}^2 H_{\sigma, \varepsilon} = 2\varepsilon + O(\sigma\varepsilon + \varepsilon^2).
\end{equation*}

Let $R$ be the cylindrical radius and $\theta$ the planar area angle. We abbreviate $H = H_{\sigma, \varepsilon}$, but keep the subscript $\varepsilon$ in the slope~$\kappa_\varepsilon$ to distinguish it from the zero-twist value~$\kappa$. Because $\zeta = -Z/N$, a slope~$\kappa_\varepsilon$ with respect to $(\Theta, \zeta)$ is a slope~$-N\kappa_\varepsilon$ with respect to $(\Theta, Z)$. We can now perform the straightening of Section~\ref{sec reference}, replacing its angular correction by $u_\theta = R^2/\langle R^2\rangle_\theta - 1$ and putting $\Theta = \theta + u$ and $Z = -\varphi - N\kappa_\varepsilon(\mathcal I) u$. The resulting formulas are
\begin{equation}
\begin{gathered}
\varrho(\mathcal I) = \sigma^2 \int_0^{\mathcal I} \langle R^2\rangle_\theta(s) \dd s, \qquad \omega = \frac{\partial_{\mathcal I} H}{\langle R^2\rangle_\theta}, \qquad \kappa_\varepsilon = \frac{\kappa}{\partial_{\mathcal I} H}, \\
b = \omega V, \qquad V = X - N\kappa_\varepsilon U, \qquad j_0 = \beta_0 U, \qquad \beta_0 = F_\sigma(H) + 8\varepsilon H.
\end{gathered}
\label{kg118}
\end{equation}
As before, passage to $\varrho$ puts the volume form in the form~$\mu = d\varrho \wedge d\Theta \wedge dZ$. We shall henceforth view $\omega, \kappa_\varepsilon, \beta_0$ as functions of $\varrho$, and a prime will mean $\partial_\varrho = (\sigma^2 \langle R^2\rangle_\theta)^{-1} \partial_{\mathcal I}$. In particular,
\begin{equation*}
\kappa_\varepsilon' = -\frac{\kappa\,\partial_{\mathcal I}^2 H}{\sigma^2 \langle R^2\rangle_\theta(\partial_{\mathcal I} H)^2} < 0.
\end{equation*}
Notice that the outer action $\varrho_b = (\sigma^2/2\pi) \int_D e^{2\sigma x}\,dx\,dy$ does not depend on $\varepsilon$. The charts therefore identify the different source domains with one common cap. Write $\Phi_{0, \varepsilon}$ for the reference immersion in these coordinates, with $\Phi_{0, 0} = \Phi_0$. If $\sigma$ is kept fixed, the profiles, coordinate charts and metrics converge in all the required finite norms as $\varepsilon\downarrow0$.

We next adapt the inverse construction on $M = \Sigma \times \T_N$ from Section~\ref{kg41}. The comparison metric remains $g_*$, with the constant shear $c_* = 2\kappa_*$. Extend the background coefficients by Lemma~\ref{kg26}, taking $N(\kappa_\varepsilon - \kappa)$ small in the finite norms which will be specified below. For weights~$n \in \Z$ and $m \in N^{-1} \Z$, set~$k(\varrho) = n - Nm\kappa_\varepsilon(\varrho)$. The band used in~\eqref{kg46} is now replaced by a fixed cone,
\begin{equation*}
\mathsf P = \sum_{n, m\ne0, \ \abs{n-\kappa_* m}<\gamma \abs{m}} \mathsf P_{nm}, \qquad \mathsf E = \Id - \av - \mathsf P, \qquad 0 < \gamma \ll 1.
\end{equation*}
The condition~$N \norm{\kappa_\varepsilon - \kappa}_\infty < \gamma/2$ ensures that $\abs{k} \geq c_\gamma(\abs{n} + \abs{m})$ on $\mathsf E$. On $\mathsf P$, on the other hand, the ratios $m/n$ and $n/m$ are both bounded. If $j = \alpha b + \beta U$, its transport frequency~$\lambda = \alpha\omega k + \beta m$ obeys $\abs{\lambda} \geq c \abs{m} \geq c/N$ throughout the fixed low-norm coefficient neighborhood. Recall also that the transport inverse for a radial angular field~$Y$, in a mode with nonzero frequency~$\nu = nY^\Theta + mY^Z$, is given by $T_Y v = v/(\ii\nu) + v^\varrho Y_\varrho/(\ii\nu)^2$, with $Y_\varrho = (Y^\Theta)'X + (Y^Z)'U$. Thus we can still use $T_b$ on $\mathsf E$ and $T_j$ on $\mathsf P$.

We keep the same augmented unknowns $f = \omega\xi^\varrho$ and $\tau = w^Z$; the modification lies in the formula that reconstructs the magnetic variation:
\begin{equation}
w = Sz = e + \frac1nX^\Sigma_{kf} + J_2 \tau, \qquad w^\varrho = \ii kf, \qquad w^\Theta = -\frac{(kf)'} n - \frac mn\tau.
\label{kg119}
\end{equation}
In particular, the derivative now contains the extra term in $(kf)' = kf' - Nm\kappa_\varepsilon'f$. We retain the variational form~\eqref{kg59}, taking $A_E = \CC^{-1} \mathcal{L}_j T_b$ and the coefficients
\begin{equation*}
\mathsf B = -\frac{\lambda}{n\omega}, \qquad \mathsf C = \frac{\beta\omega'}{\omega^2} + \frac{Nm}{n} \frac{\beta\kappa_\varepsilon'}{\omega} + \frac kn\left[\left(\frac\beta\omega\right)' - N\kappa_\varepsilon \alpha'\right].
\end{equation*}
Testing against the scalar components gives $(k/n) \mathfrak r^Z + \mathsf Cf + \mathsf B\tau = 0$ and $-\ii\mathfrak r^\varrho/n + \mathsf Bf = 0$, where we have denoted $\mathfrak r = (\Id - \av) \allowbreak\curl_\mu(G(w, \cdot)) - g$. Substituting these identities in the reconstruction shows that $\xi = T_b e + T_j \mathsf P\mathfrak r$ satisfies both $[b, \xi] = w$ and $[j, \xi] = \mathfrak r$. Notice that the computation is valid also at $k(\varrho) = 0$, so no rational pressure surface has been removed.

The first issue is coercivity, since more than one poloidal weight may now occur for a fixed longitudinal weight~$m$. Choose a smooth real function~$\chi = 1$ on $[-1, 1]$, supported in $(-2, 2)$ and satisfying $0 \leq \chi \leq 1$. Put $\mathsf T = \mathsf P\chi(k/\delta)$, with $0 < \delta < \min(1, \gamma, \kappa_* - \gamma/2)/(8N)$. For a tangent field~$v_{n'}$ of poloidal weight~$n'$, consider the Hodge decomposition $v_{n'} = \Pi_\Sigma v_{n'} + \nabla_\Sigma p_{n'}$, where $p_{n'}$ has zero mean. The spectral gap on the sphere implies
\begin{equation*}
\norm{p_{n'}}_0 \leq \frac{CR_\Sigma}{1 + \abs{n'}} \norm{v_{n'}}_0, \qquad \norm{\nabla_\Sigma \Delta_\Sigma^{-1} \phi_{n'}}_0 \leq \frac{CR_\Sigma}{1 + \abs{n'}} \norm{\phi_{n'}}_0
\end{equation*}
whenever $\phi_{n'}$ is a scalar of the same weight; in the case $n' = 0$, its constant mode is omitted. For a radial multiplier $\chi_n$ we have the commutator identity
\begin{equation*}
[\Pi_\Sigma, \chi_n] v_{n'} = -\Pi_\Sigma(p_{n'} \nabla_\Sigma \chi_n) - (\Id - \Pi_\Sigma)(\chi_n \Pi_\Sigma v_{n'}).
\end{equation*}
The potential corresponding to the second projection has Laplacian $-(\Pi_\Sigma v_{n'}) \cdot \nabla_\Sigma \chi_n$. This has zero mean, including in the case $n' = 0$, so the same spectral estimate gives
\begin{equation}
\norm{[\Pi_\Sigma, \chi_n] v_{n'}}_0 \leq \frac{CR_\Sigma}{1 + \abs{n'}} \norm{\nabla_\Sigma \chi_n}_\infty \norm{v_{n'}}_0.
\label{kg120}
\end{equation}

Fix a longitudinal weight~$m$, and use the notation $\chi_n = \chi((n - Nm\kappa_\varepsilon)/\delta)$ and $\eta_\kappa = NR_\Sigma^2 \norm{\kappa_\varepsilon'}_\infty/\delta$. Let $\phi = \sum_n \phi_n$ be a scalar in the cone belonging to this longitudinal sector. Combining the bounds $\norm{\nabla_\Sigma \chi_n}_\infty \leq CNR_\Sigma \abs{m} \norm{\kappa_\varepsilon'}_\infty/\delta$ and $\abs{m}/(1 + \abs{n \pm 1}) \leq C_{\kappa_*}$ with~\eqref{kg120}, we obtain
\begin{equation*}
\norm{[\Pi_\Sigma, \chi_n](\phi_n \nabla_\Sigma h_\Sigma)}_0 \leq C_{\kappa_*} \eta_\kappa \norm{\phi_n}_0.
\end{equation*}
Define $u_{\Sigma, n} = \Pi_\Sigma(\chi_n \phi_n \nabla_\Sigma h_\Sigma)$. By Lemma~\ref{kg65}, we have $\sum_n \norm{u_{\Sigma, n}}_0^2 \leq \frac12 \norm{\phi}_0^2$. Now $u_{\Sigma, n}$ contains only the weights~$n \pm 1$, which means that a cross term can arise only between $n$ and $n \pm 2$. For these pairs we have $\chi_n \chi_{n'} = 0$, because the two frequencies differ by two whereas $4\delta < 2$. It follows that
\begin{equation*}
\langle u_{\Sigma, n}, u_{\Sigma, n'}\rangle = \langle[\Pi_\Sigma, \chi_n](\phi_n \nabla_\Sigma h_\Sigma), \chi_{n'} \phi_{n'} \nabla_\Sigma h_\Sigma\rangle.
\end{equation*}
There are only two possible neighbours of a given $n$. Summing the preceding cross terms therefore gives at most $C_{\kappa_*} \eta_\kappa \norm{\phi}_0^2$, with no dependence on the number of selected weights. Consequently, if we write the Hodge decomposition as $(\mathsf T\phi) \nabla_\Sigma h_\Sigma = u_\Sigma + \nabla_\Sigma p_\Sigma$, normalizing $p_\Sigma$ to have zero mean, then
\begin{equation}
\norm{u_\Sigma}_0^2 \leq (\tfrac12 + C_{\kappa_*} \eta_\kappa) \norm{\phi}_0^2, \qquad \norm{mp_\Sigma}_0 \leq C_{\kappa_*} R_\Sigma \norm{\phi}_0.
\label{kg121}
\end{equation}
For the second inequality, apply the Hodge decomposition separately to each field~$\chi_n \phi_n \nabla_\Sigma h_\Sigma$. Its potential $p_{\Sigma, n}$ satisfies $\norm{p_{\Sigma, n}}_0 \leq CR_\Sigma \abs{n}^{-1} \norm{\chi_n \phi_n}_0$ by Lemma~\ref{kg65}, and a given output weight receives contributions from at most two such fields. Summing first over finitely many weights and then passing to the limit in $L^2$ proves the assertion.

Consider a solenoidal field~$w = v + w^Z U$ in this sector and take $\phi = \mathsf Tw^Z$. We use the energy $\mathcal E_m = \norm{v}_0^2 + \norm{w^Z + c_* dh_\Sigma(v)}_0^2$, as in the proof of Lemma~\ref{kg67}. The self-adjointness of $\mathsf T$ then gives
\begin{equation*}
\norm{\phi}_0^2 = \langle w^Z, \mathsf T\phi\rangle = \langle w^Z + c_* dh_\Sigma(v), \mathsf T\phi\rangle - c_* \langle v, u_\Sigma\rangle - c_* \langle v, \nabla_\Sigma p_\Sigma\rangle.
\end{equation*}
Estimate~\eqref{kg121} bounds the first two terms by $\mathcal E_m^{1/2}(1 + c_*^2/2 + C_{\kappa_*} \eta_\kappa)^{1/2} \norm{\phi}_0$. In the last term, integration by parts together with $\Div_\Sigma v = -\ii mw^Z$ and $\norm{w^Z}_0 \leq (1 + c_*^2)^{1/2} \mathcal E_m^{1/2}$ gives $C_{\kappa_*} R_\Sigma \mathcal E_m^{1/2} \norm{\phi}_0$. We divide by the scalar norm, square, and sum orthogonally over $m$. For $R_\Sigma, \eta_\kappa \leq 1$, the result is
\begin{equation}
\norm{\mathsf Tw^Z}_0^2 \leq \left(1 + 2\kappa_*^2 + C_{\kappa_*} \left(R_\Sigma + \frac{NR_\Sigma^2 \norm{\kappa_\varepsilon'}_\infty}{\delta}\right)\right) \int_M g_*(w, \overline w)\,\mu.
\label{kg122}
\end{equation}
On the selected scalar components, set~$\mathsf W = (\Id - \mathsf T^2)^{1/2}$. The relations $(\mathsf Pw)^\varrho = \ii kf$ and $\chi_n = 1$ for $\abs{k} \leq \delta$ give
\begin{equation*}
\norm{\mathsf Wf}_0 \leq (R_\Sigma/\delta) \norm{w}_0, \qquad \norm{\mathsf W\tau}_0 \leq \norm{w}_0.
\end{equation*}
We now split the scalar part of the variational form according to $\mathsf T$ and $\mathsf W$. In the first part, completion of the square yields
\begin{equation*}
\norm{\mathsf C^{1/2}(\mathsf Tf + \mathsf C^{-1} \mathsf B\mathsf T\tau)}_0^2 - \langle\mathsf B^2 \mathsf C^{-1} \mathsf T\tau, \mathsf T\tau\rangle,
\end{equation*}
and the second part has the lower bound $-C(R_\Sigma/\delta + R_\Sigma^2/\delta^2) \norm{w}_0^2$. On the support of $\mathsf T$ the coefficient~$\mathsf C$ is positive, with $\mathsf B^2/\mathsf C = 1/(3\kappa_*^2) + o(1)$. Combining this fact with~\eqref{kg122}, the metric comparison and $\norm{A_E} \leq C_{\gamma, N} R_\Sigma$, we recover the positive margin $(\kappa_*^2 - 1)/(3\kappa_*^2)$, up to small errors. The parameters are chosen in the following order: first $L$ large, then $\gamma, \delta$ small, then $\sigma$ small, and finally the twist and the coefficient neighborhood sufficiently small. The square left above controls $\mathsf Tf$ by the estimate for $\mathsf T\tau$, while the displayed bound controls $\mathsf Wf$. Using~\eqref{kg55}, we conclude that the form is coercive on the completion of smooth triples for the norm~$E_0$.

We must now prove the finite regularity estimates on this same low-norm neighborhood. Fix a finite derivative reserve $K$, and on the cone write $H_\Sigma f = X_f^\Sigma$ and $M_1 = m/n$. We use the product grading~$\Lambda^s$ introduced in Section~\ref{kg41}. The identity~$J_1 = H_\Sigma(\Id - N\kappa_\varepsilon M_1)$ shows that the commutator splits as
\begin{equation*}
\mathcal K_s z = (\Lambda_1^s S - S\Lambda^s) z = \mathcal K_s^{(1)} f + \mathcal K_s^{(2)} \tau, \qquad \mathcal K_s^{(1)} = -NH_\Sigma[\Lambda_0^s, \kappa_\varepsilon] M_1.
\end{equation*}
Here $\mathcal K_s^{(2)}$ is the contribution from $J_2$. At an even order~$2 \leq s \leq K$, the grading~$\Lambda_0^s = (1 + \Delta_0)^{s/2}$ is differential, so $\mathcal K_s^{(1)}$ and its adjoint $-NM_1[\kappa_\varepsilon, \Lambda_0^s] H_\Sigma^*$ have order at most $s$; their coefficients contain derivatives of $\kappa_\varepsilon$. In particular, if $a_K = NC_K \norm{\kappa_\varepsilon - \kappa}_{C^{K+2}}$, we have
\begin{align}
\bigl\|\mathcal K_s^{(1)} f\bigr\|_0 &\leq a_K \norm{f}_s, & \bigl\|(\mathcal K_s^{(1)})^* v\bigr\|_0 &\leq a_K \norm{v}_s, \notag \\
\bigl\|\mathcal K_s^{(2)} \tau\bigr\|_0 &\leq C_s \norm{\tau}_{s-1}, & \bigl\|(\mathcal K_s^{(2)})^* v\bigr\|_0 &\leq C_s \norm{v}_{s-1}.
\label{kg123}
\end{align}
All the norms in these estimates are formed with the fixed product metric; their constants may depend on the thickness. Consider the graph norms $E_s$ from~\eqref{kg73}. The identities $e = \mathsf Ew$ and $\tau = \mathsf P(w^Z)$ allow us to repeat the argument of~\eqref{kg55}, obtaining $E_s(z)^2 \sim \norm{w}_s^2 + \norm{f}_s^2$ for $w = Sz$. In particular, $\norm{\Lambda^s z}_0 \leq E_s(z)$. The contribution of $J_1 f$ is the only one that differs between $E_0(\Lambda^s z)$ and $E_s(z)$. The estimate for $\mathcal K_s^{(1)}$ in~\eqref{kg123} therefore shows that, when $a_K$ is sufficiently small,
\begin{equation}
E_0(\Lambda^s z)^2 \geq cE_s(z)^2, \qquad \norm{S\Lambda^s z}_0 \leq CE_s(z),
\label{kg124}
\end{equation}
where $c > 0$ is fixed.

We shall keep track of the high coefficient norms on the right-hand side. For the coefficients~$p = (G, \alpha, \beta)$, choose an even low index~$s_0$ sufficiently large and a fixed shift $d_0$, both independent of $K$, and restrict the coefficient neighborhood by $\norm{p}_{s_0+d_0} \leq M_0$. As in Section~\ref{kg41}, this norm includes the radial derivatives appearing in the form. We identify $G$ with its endomorphism relative to the product metric. The new adjoint contribution in~\eqref{kg81} then satisfies
\begin{equation}
\bigl\|(\mathcal K_s^{(1)})^*(Gw)\bigr\|_0 \leq C_s a_K\bigl(M_0 E_s(z) + \norm{G}_s E_{s_0}(z)\bigr),
\label{kg125}
\end{equation}
by the Sobolev product estimate. Thus the high norm of $G$ is multiplied by $E_{s_0}(z) E_s(z)$ in the energy identity, while the other new contribution satisfies $\bigl|\langle G\mathcal K_s^{(1)} f, S\Lambda^s z\rangle\bigr| \leq C_s a_K M_0 E_s(z)^2$.

For all the other commutators, the bounds used in~\eqref{kg81} remain available. In the off-cone block, we factor $n - \kappa_* m$ out of $k$ and expand the reciprocal in powers of $Nm(\kappa_\varepsilon - \kappa)/(n - \kappa_* m)$. Since this ratio has modulus smaller than $1/2$, and differentiation introduces only polynomial factors in the summation index, the differentiated series are summable through the chosen finite reserve. Applying the product and commutator estimates of Section~\ref{kg41}, we obtain, for $s \geq s_0$,
\[\norm{[\Lambda_1^s, G] w}_0 + \bigl\|(\mathcal K_s^{(2)})^*(Gw)\bigr\|_0 + \norm{[\Lambda^s, \mathsf R_p] z}_0 \leq C_s\bigl(E_{s-1}(z) + (1 + \norm{p}_{s+d_0}) E_{s_0}(z)\bigr),\]
where $\mathsf R_p$ denotes the collection of blocks $-A_E, \mathsf B, \mathsf C$. Substitution of these bounds, together with~\eqref{kg125}, in the energy identity~\eqref{kg81} gives
\[\bigl|\cform_p(z, \Lambda^{2s} z) - \cform_p(\Lambda^s z, \Lambda^s z)\bigr| \leq C_{K, M_0} a_K E_s(z)^2 +C_s\bigl(E_{s-1}(z) + (1 + \norm{p}_{s+d_0}) E_{s_0}(z)\bigr) E_s(z).\]
We have used $a_K \leq 1$ in the second term. We can now choose the twist small enough that the first term is absorbed by less than half of the positive constant coming from coercivity and~\eqref{kg124}. This choice depends on $K$ and the bound on the low coefficient norm, and is independent of the high norms of $G$.

For the weak equation~\eqref{kg74}, retain the data norm~$D_s = \norm{H_1}_s + \norm{H_2}_s$. Taking $\Lambda^{2s} z$ as test function gives the right-hand side
\begin{equation*}
\langle\Lambda_1^s H_1, S\Lambda^s z\rangle - \langle\mathcal K_s^* H_1, \Lambda^s z\rangle + \langle\Lambda^s H_2, \Lambda^s z\rangle,
\end{equation*}
which is bounded in absolute value by $C_s D_s E_s(z)$, in view of~\eqref{kg123}--\eqref{kg124}. We have therefore proved
\begin{equation*}
E_s(z) \leq C_s\bigl(D_s + E_{s-1}(z) + (1 + \norm{p}_{s+d_0}) E_{s_0}(z)\bigr).
\end{equation*}
For the finitely many even orders not exceeding $s_0$, the bounded low coefficient derivatives give the sharper lower-order contribution $C_s E_{s-1}(z) E_s(z)$, once the same absorption has been made. Thus the estimate at those orders reduces to $E_s(z) \leq C_s(D_s + E_{s-1}(z))$. Applying interpolation directly to the pair~$(w, f)$, we have
\begin{equation*}
E_{s-1}(z) \leq \nu E_s(z) + C_{s, \nu} E_0(z).
\end{equation*}
Combining this with the coercive bound $E_0(z) \leq CD_0$, we first obtain $E_{s_0}(z) \leq CD_{s_0}$ and then, at every even order~$s \geq s_0$ for which $s + d_0 \leq K$,
\begin{equation}
E_s(z) \leq C_s\bigl(D_s + (1 + \norm{p}_{s+d_0}) D_{s_0}\bigr).
\label{kg127}
\end{equation}
At an arbitrary order we use the next even one, so at most two additional derivatives are lost. This loss is independent of $K$.

To justify these calculations for weak solutions, we work first on finite spectral spaces of triples. These spaces preserve the cone and the test $\Lambda^{2s} z$, so none of the estimates depends on the Galerkin dimension. Their union is dense in the graph completion: smooth triples are dense by its definition, and their spectral approximations converge in $H^1$, hence also in the graph norm~$E_0$. Weak compactness then gives a limit, and uniqueness identifies it with the coercive solution. This proves~\eqref{kg127} for that solution.

We can consequently carry out the physical reconstruction from Section~\ref{kg41}. It gives $M_{\Phi, \varepsilon} R_{\Phi, \varepsilon} = \Id$, with fixed integers $s_b, d$ that are independent of $K$, and the estimate
\begin{equation}
\norm{R_{\Phi, \varepsilon} h}_s \leq C_s\left(\norm{h}_{s+d} + (1 + \norm{\Phi - \Phi_{0, \varepsilon}}_{s+d}) \norm{h}_{s_b}\right), \qquad \max(s + d, s_b) \leq K.
\label{kg128}
\end{equation}
On the chosen finite range of orders, the constants remain uniform as $\varepsilon\downarrow0$. The divergence and potential homotopies have the same derivative losses as before. For the transport factors, finite product estimates follow from the preceding series for $1/k$ off the cone and the identity $\lambda^{-1} = m^{-1}(\beta + \alpha\omega k/m)^{-1}$ on the cone. The reciprocal in the latter expression also expands in fixed bounded angular multipliers, and the corresponding estimates involve only one high coefficient norm. Finally, the factors $1/m$ are bounded by $N$ on the cone.

It remains to choose the derivative reserve for the nonlinear argument. Take $s_* > r + 6$ and choose an integer~$\ell$ large enough to include all the fixed losses and low indices used above, as well as those in~\eqref{kg102}. In the shifted grading $\abs{u}_s = \norm{u}_{H^{s_*+\ell s}}$, fix $K$ so that it covers every estimate through level $80$ and the further $d + 2$ derivatives required for the kernel correction. Having fixed this reserve, choose the thickness and an untwisted kernel direction~$\xi_0$ as in Section~\ref{kg90}, for the braid already selected. On the common source domain, define
\begin{equation*}
\xi_\varepsilon = \xi_0 - R_{\Phi_{0, \varepsilon}, \varepsilon} M_{\Phi_{0, \varepsilon}, \varepsilon} \xi_0.
\end{equation*}
By construction, $M_{\Phi_{0, \varepsilon}, \varepsilon} \xi_\varepsilon = 0$. The defect converges to zero in each of the required norms and involves two derivatives of $\xi_0$. Estimate~\eqref{kg128} therefore gives $\xi_\varepsilon \to \xi_0$ through the order needed to initialize the iteration. In particular, the normal trace on the reference axis still represents the chosen braid. We now fix this positive twist for the remainder of the argument.

We integrate this direction by a finite iteration. Let $\mathcal F_\varepsilon$ denote the map~\eqref{eq:nonlinear-map} associated with the reference equilibrium with twist, and introduce $\Phi = \Phi_{0, \varepsilon} + u$, $F(u) = \mathcal F_\varepsilon(\Phi)$ and $Q_u h = D\Phi\,R_{\Phi, \varepsilon} h$. In the shifted grading, the estimates take the form
\begin{equation*}
\abs{F(u)}_s \lesssim_s1 + \abs{u}_{s+1}, \qquad \abs{Q_u h}_s \lesssim_s \abs{h}_{s+1} + (1 + \abs{u}_{s+1}) \abs{h}_1.
\end{equation*}
On a fixed sufficiently small~$\abs{\cdot}_1$ ball, the derivative and the Taylor remainder are bounded in norm zero by $C \abs{v}_1$ and $C \abs{v}_1^2$, respectively. The approximate-inverse identity~\eqref{kg102} yields $\abs{DF(u) Q_u h - h}_0 \leq C \abs{F(u)}_1 \abs{h}_1$. To construct the smoothing operators, extend from the source by Lemma~\ref{kg26}, apply a spectral cutoff at scale $\mathcal N^{1/\ell}$, and restrict back; denote the result by $\mathcal S_{\mathcal N}$. A comparison of Sobolev weights on the two sides of the cutoff gives a bound $\mathcal N^{s-a}$ from level $a$ to $s \geq a$, and the tail estimate~$\mathcal N^{a-s}$ in the opposite direction. With initial value~$u_0 = tD\Phi_{0, \varepsilon} \xi_\varepsilon$, we iterate according to
\begin{equation*}
\mathcal N_{j+1} = \mathcal N_j^{3/2}, \qquad u_{j+1} = u_j - \mathcal S_{\mathcal N_j} Q_{u_j} F(u_j).
\end{equation*}
We close the induction with the bounds $\abs{F(u_j)}_0 \leq \mathcal N_j^{-16}$ and $\abs{u_j}_{80} \leq \mathcal N_j^{44}$, which imply $\abs{F(u_j)}_1 \lesssim \mathcal N_j^{-15}$ by interpolation with level $79$. Thus the increments are bounded by $C_s \mathcal N_j^{s-15}$. The inverse defect, the smoothing tail estimated at level $78$, and the Taylor remainder are bounded, respectively, by $C\mathcal N_j^{-30}$, $C\mathcal N_j^{-33}$ and $C\mathcal N_j^{-28}$. These three estimates imply the next residual bound $\mathcal N_{j+1}^{-16} = \mathcal N_j^{-24}$. The high-norm induction also closes, because $44 < 65 < 66 = (3/2) 44$. At level one the increments form a summable series and keep the iterates inside the fixed low-norm ball. The derivative estimate makes $F$ continuous from level one to level zero, so their limit satisfies $F = 0$. Finally, since $F(u_0) = O(t^2)$, we can choose $\mathcal N_0 = c \abs{t}^{-1/8}$, with $c > 0$ sufficiently small, to obtain
\begin{equation*}
F(u_\infty) = 0, \qquad \Phi_{t, \varepsilon} = \Phi_{0, \varepsilon} + tD\Phi_{0, \varepsilon} \xi_\varepsilon + O(\abs{t}^{15/8}) \quad \hbox{in }H^{s_*}.
\end{equation*}
Upon division by $t$, the normal trace converges in $C^1$ to the prescribed braid profile. The argument of Section~\ref{kg106} therefore produces an embedded axis and an embedded pressure sublevel. Both the pressure formula~\eqref{kg34} and the force identity~\eqref{kg35} are valid with $b = \omega(X - N\kappa_\varepsilon U)$. Our choice of derivative reserve gives $B, P \in C^r$ through the axis, together with a boundary of class $C^{r+1}$. In the angular coordinate~$\zeta = -Z/N$, the flow has slope~$\kappa_\varepsilon$; formula~\eqref{kg118} gives $\kappa_\varepsilon \to \kappa$ as $\varepsilon\downarrow0$, and the derivative calculated above satisfies $\kappa_\varepsilon' < 0$.
\end{proof}

The twist must be small in norms involving derivatives up to $K + 2$. Since the required smallness changes when the desired regularity increases, the conclusion is a nonzero twist for each prescribed finite regularity. The argument does not provide a $C^\infty$ construction with a single fixed nonzero twist.

\section*{Acknowledgements}

This work has received funding from the European Research Council (ERC) under the European Union's Horizon 2020 research and innovation programme through grant agreement 862342 (A.E.). A.E. and D.P.-S. also acknowledge support from the grants \nolinkurl{PID2022-136795NB-I00}, \nolinkurl{PID2025-168410NB-I00} and \nolinkurl{CEX2023-001347-S} (ICMAT--Severo Ochoa), funded by \nolinkurl{MCIN/AEI/10.13039/501100011033}. J.P.-T. was supported by the Swiss State Secretariat for Education, Research and Innovation (SERI) under contract number MB22.00034 through the project TENSE.

\section*{Declaration on the use of AI tools}

ChatGPT was used during the preparation of this work to assist with mathematical exploration and calculations, bibliographic searches, and the organization, drafting and revision of the manuscript. The authors take full responsibility for the mathematical results, arguments and final content of this paper.

\end{document}